\documentclass[10pt]{article}%
\usepackage{amsfonts,amssymb}
\usepackage{amsmath}
\usepackage{graphicx,epsfig}
\usepackage{color}
\usepackage{multirow}
\usepackage{hyperref}
\usepackage{indentfirst}
\newtheorem{thm}{Theorem}[section]

\newtheorem{lem}[thm]{Lemma}

\newtheorem{defn}[thm]{Definition}
\newtheorem{rem}[thm]{Remark}
\newtheorem{exm}[thm]{Example}

\numberwithin{equation}{section}

 \newenvironment{proof*}[1]{\par \textbf{Proof.}\quad #1}{\hfill \textbf{} $\Box$ }

\newcommand{\px}[1][x]{\partial_{#1}}
\title{Implicit-Explicit difference schemes for nonlinear fractional differential equations with non-smooth solutions}

\author{Wanrong Cao\footnotemark[1], Fanhai Zeng\footnotemark[2], Zhongqiang Zhang\footnotemark[3]
,  and George Em Karniadakis \footnotemark[4]}
\begin{document}
\maketitle
\renewcommand{\thefootnote}{\fnsymbol{footnote}}

\footnotetext[1]{Department of Mathematics, Southeast University, Nanjing 210096, P.R.China. (wrcao@seu.edu.cn)}
\footnotetext[2]{ Division of Applied Mathematics, Brown University, Providence
RI, 02912, USA. (fanhai\_zeng@brown.edu)}
\footnotetext[3]{Department of Mathematical Sciences, Worcester Polytechnic Institute, Worcester, MA, 01609, USA.  (zzhang7@wpi.edu).}
\footnotetext[4]{ Division of Applied Mathematics, Brown University, Providence
	RI, 02912, USA.  (george\_karniadakis@brown.edu).}

\begin{abstract}
We propose second-order implicit-explicit (IMEX) time-stepping schemes  for nonlinear fractional differential equations with fractional order $0<\beta<1$. From the known structure of the non-smooth solution and by introducing corresponding correction terms,
we  can obtain uniformly second-order accuracy from  these schemes. We   prove the convergence and linear stability  of the proposed schemes.
 Numerical examples illustrate the flexibility and efficiency of the IMEX schemes and  show that they are effective for  nonlinear and multi-rate fractional  differential systems as well as multi-term fractional differential systems with non-smooth solutions.

\vspace{0.5em} 

\end{abstract}

%%%%%%%%%%%%%%%%%%%%%%%%%%%%%%%%%%%%%%%%%%%%%
%%%%%%%%%%%%%%%%%%%%%%%%%%%%%%%%%%%%%%%%%%%%%

\section{Introduction}

%%%%%%%%%%%%%%%%%%%%%%%%%%%%%%%%%%%%%%%%%%%%%
%%

We aim at constructing efficient  finite difference schemes for fractional ordinary differential equations (FODEs) with non-smooth solutions.
In recent decades, due to the increasing interest in problems with anomalous transport dynamics, fractional differential equations have become significant mathematical models in many fields of science and engineering, such as viscoelastic models in blood flow \cite{ParisGeorge2014}, underground transport \cite{Hatano98}, options pricing model in financial markets \cite{WangChen2015}, etc.

Though some fractional differential equations (FDEs) with special form, e.g., linear equations, can be solved by analytical methods, e.g., the Fourier transform method or the Laplace transform method \cite{Podlubny1999}, the analytical solutions of many generalized FDEs (e.g. nonlinear FDEs and  multi-term FDEs) are rather difficult to obtain. This encourages us to develop effective numerical methods for solutions of these FDEs.
Up to now, a number of finite difference methods have been established for FDEs. One way is to transform the considered  FDEs into their integral forms, then numerical methods for the fractional integral operator are developed and the corresponding difference schemes are derived; see    \cite{CuestaLub2006,CuePal03,DieFF02,DieFF04,Gar10,Garrappa2015,GalGar08,Lub85,Lub86,Tang93,ZengLLT13}.
Another approach is to approximate the fractional derivative operators in the considered FDEs  directly; see   \cite{GaoSZ14,LinXu07,SunWu06,TianDeng2015,ZhouDeng2013}.  Besides finite difference methods, there exist also  other numerical  methods for FDEs, e.g., finite element methods \cite{JinLPZ2014,WangYang2013}, spectral methods \cite{LiXu2009,TianDW2014,MoshenGeorge2014}, matrix methods \cite{PanRM2014,PodlubnyCS2009}, etc.

For nonlinear FDEs, most of the aforementioned finite difference methods are implicit, and the nonlinear system needs to be solved using the iteration method, which is costly.
To avoid extra computational costs and instability caused by iteration, various numerical methods have been proposed, e.g., predictor-corrector methods; see   \cite{DafYS14,Deng07,DieFF02, DieFF04,Gar10,LiCY11,YangLiu05,ZhaoDeng14}  and time-splitting schemes \cite{CaoZhangK15}.

Implicit-explicit (IMEX) schemes (also known as semi-implicit schemes, linearly-implicit  schemes) play a major rule in the numerical treatment of nonlinear/stiff differential equations \cite{Giacomoimex2013,Heimex2012,Knio1999}.  Due to approximating the nonlinear/nonstiff part explicitly, we avoid solving nonlinear equations/systems in every time step, and hence the computational cost can be reduced significantly. Moreover, the IMEX schemes generally have better stability than explicit schemes, for the linear/stiff parts are treated implicitly. To the best of authors' knowledge, there are very limited works on IMEX schemes for FDEs.
An IMEX scheme has been proposed for nonlinear fractional anomalous diffusion equations with smooth solutions in \cite{ZengLL13} and an explicit and implicit finite difference scheme for fractional Cattaneo equation was developed in \cite{GhaMae2010}. Some semi-implicit methods for space-fractional differential equations can be found in the literature, e.g., see \cite{ChenLiu-etal15,LiuZhuang15,Yang16}.

The main contribution of this work is to develop second-order IMEX schemes for nonlinear/stiff FODEs with solutions that have a weak singularity at the origin. Most numerical methods for differential equations are generally intended for problems with solutions of high regularity. However, solutions to FDEs usually have weak singularities at the origin even when the forcing term is smooth, see e.g.  \cite{DieFor02,DieFF04,Kilbas2006}. When solving these FDEs, the singularity requires special attention to obtain the expected high accuracy. Several approaches have been proposed to deal with the weak singularity,  such as using adaptive grids (nonuniform grids) to keep errors small near the singularity \cite{Garrappa2015,YusteJoaquin12,ZhangSunLiao2014},  or employing non-polynomial basis functions to include the correct singularity index \cite{CaoHX2003,FordMR2013,ZhangZK2015}, or using the correction terms to remedy the loss of accuracy and  recover high-order schemes \cite{Garrappa2015,Lub86,Zeng2015,ZengLLT2015,ZengZK2015}.

In this paper, we follow the last approach
to develop IMEX methods with uniformly second-order accuracy for FODEs, whose solutions are non-smooth and have known structure. To deal with the singularity,  we apply correction terms in the proposed schemes, so that the resulting fractional quadratures are either exact or sufficiently accurate for the weakly singular parts of the solutions. The idea of adding correction terms was firstly proposed for approximating fractional-order integrals by the linear multi-step method   in \cite{Lub86}. Very recently, the same strategy has been adopted to enhance the accuracy of numerical schemes for FDEs  \cite{ZengLLT2015,ZengZK2015}. To treat the nonlinear part explicitly, we further use extrapolation and Taylor expansion to approximate the nonlinear part, where appropriate correction terms are also used to obtain high accuracy for non-smooth solutions.
We propose two IMEX schemes, which can work for nonlinear/stiff FODEs with uniformly second-order convergence, even when the solutions have weak singularity at the origin. We also prove the convergence and linear stability of the proposed schemes.

We organize this work as follows. In Section 2, we formulate the IMEX method based on the  fractional linear multistep method (FLMM), and then derive two IMEX schemes by applying the extrapolation and Taylor expansion for the nonlinear terms. Moreover, we provide the strategies for introducing correction terms to the schemes to make them uniformly second-order for FODEs with non-smooth solutions. We also present convergence rates of the proposed schemes, the proofs of which are given in Section 6.  In Section 3, we discuss the linear stability of the proposed schemes. We present more details on the proposed IMEX schemes in Section 4 to show that they can be applied to stiff systems and multi-term nonlinear FODEs.  In Section 5, we present numerical examples to illustrate the computational flexibility and verify our error estimates. We conclude in Section 7 and discuss the performance of the proposed schemes.

%%%%%%%%%%%%%%%%%%%%%%%%%%%%%%%%%%%%%%%%%%%
\section{Second-order IMEX schemes}\label{sec:time-splitting-schemes}
We consider the following nonlinear FODE
\begin{eqnarray}\label{eq:fde}
(\,^CD_{0}^{\beta}u)(t)&=& \lambda u(t)+f_u(t),\; t\in(0,T],\; u(0)=u_0,
\end{eqnarray}
where $0<\beta<1$,
%$A$ is a $d\times d$ real-valued matrix,
%and $f:[0,T]\times \mathbb{R}^d\to\mathbb{R}^d$,
$f_u(t)=f(t,u(t))$,  ${}^CD_{0}^{\beta}$  is the Caputo derivative defined by
\begin{equation}\label{eq:caputo1}
(\,^CD_{a}^{\beta} g)(t) =(I_{a}^{1-\beta} g')(t), \quad (I_{a}^\beta g)(t)=\frac{1}{\Gamma(\beta)}
\int_a^t\frac{g(\tau)}{(t-\tau)^{1-\beta}}\mathrm{d}\tau,~t>a.
\end{equation}

We  first transform  \eqref{eq:fde} into its integral form as
\begin{eqnarray}\label{eq:int}
{u}(t)&=&{u}(0)+\lambda\,(I_{0}^\beta u)(t) +(I_{0}^\beta f_u)(t),\quad 0\leq t\leq T.
\end{eqnarray}
Eq. \eqref{eq:int} is readily obtained by applying the operator $I_{0}^\beta$ on both sides
of \eqref{eq:fde} and using the identity $(I_{a}^\beta \,^CD_{a}^\beta g)(t)=g(t)-g(a)$,
see, e.g. \cite{DieFF04}.

In order to obtain second-order schemes, we need to approximate the
fractional integrals in \eqref{eq:int} with second-order quadrature rules. However, the solutions to \eqref{eq:fde} usually  have singularity at $t=0$.  The analytical solution to \eqref{eq:fde} can be written as the summation of  regular  and  singular parts, as given in the following lemma. See also \cite{DieFFW06,DieFor02,DieFF04} for more discussions.
\begin{lem}[\cite{DieFF04}]\label{lem:diethm1}
Suppose $f\in\mathcal{C}^3(G)$, where  $G$ is a suitable region of variable $u$. Then there exists a function $\psi\in\mathcal{C}^2[0,T]$ and some $c_1,\cdots,c_{\hat{\nu}}\in \mathbb{R}$ and $d_1,\cdots,d_{\hat{\eta}}\in \mathbb{R}$, such that the solution of \eqref{eq:fde} is of the form
\begin{equation}\label{eq:structure}	
u(t)=\psi(t)+\sum_{\nu=1}^{\hat{\nu}}c_{\nu}t^{\nu\beta}
+\sum_{\eta=1}^{\hat{\eta}}d_{\eta}t^{1+\eta\beta},
\end{equation}
where $\hat{\nu}:=[2/\beta]-1$, $\hat{\eta}:=[1/\beta]-1$.
\end{lem}

The solution of   \eqref{eq:fde} is usually non-smooth, even if $f(t,u(t))$
is smooth. Consequently, many existing numerical methods (see e.g.\cite{DieFF02,Garrappa2015,ZengLL13}) for  \eqref{eq:fde} would produce less accurate numerical solutions when they are directly applied.
Next, we will adopt a second-order  FLMM  developed in \cite{Lub86} to \eqref{eq:int} to construct our IMEX schemes. Take a uniform partition of time interval $[0,T]$, i.e., $t_n=nh,0\leq n\leq N$ with $h=T/N$.
The second-order  FLMM used in the present work reads
\begin{equation}\label{eq:lubapp1}
I_{0}^{\beta} u(t_n)=h^\beta \sum_{j=0}^n\omega_{n-j}^{(\beta)} u(t_j)
+{h^\beta} \sum_{j=0}^mW_{n,j}^{(\beta)} u(t_j) + O(h^2),
\end{equation}
where $\{\omega_{j}^{(\beta)}\}$ are coefficients of the Taylor expansion of the following  generating function
\begin{equation}\label{eq:genfun}
\omega^{(\beta)}(z)=\left(\frac{1}{2}\frac{1+z}{1-z}\right)^\beta=\sum_{j=0}^\infty
\omega_{j}^{(\beta)}z^j,
\end{equation}
and $\{W_{n,j}^{(\beta)} \}$ are the starting weights that  recover  second-order accuracy.  If we drop the correction terms $\sum_{j=0}^mW_{n,j}^{(\beta)} u(t_j)$
in \eqref{eq:lubapp1}, then we would lose the second-order accuracy, unless $u(t)$ satisfies some
special conditions, i.e., $u(t)$ is smooth and $u(0)=u'(0)=0$.

The following lemma states the convergence of \eqref{eq:lubapp1} with no correction terms.
\begin{lem}[\cite{Lub86,ZengLLT13}]\label{lem:lubich}
If $u(t)=t^{\nu},\;\nu\geq 0$, then for \eqref{eq:lubapp1} with $m=0$,
\begin{eqnarray}\label{eq:lubich-1}
(I_0^\beta u)(t_n)=h^\beta \sum_{k=0}^n\omega_{n-k}^{(\beta)} u(t_k)
+O(h^2t_n^{\nu+\beta-2})+O(h^{1+\nu}t_n^{\beta-1}),
\end{eqnarray}
where $\{\omega_{k}^{(\beta)}\}$ are defined by \eqref{eq:genfun}.	
\end{lem}

%By the lemma, $(\mathcal{I}_h^\beta g)$ is of second-order accuracy uniformly for $0<n\leq N$ only when  $\nu\geq 2-\beta$. For small $\nu\geq 0$,   the approximation has only second-order accuracy when $t_n$ is not close to $t_0$.

From Lemma \ref{lem:diethm1}, we see that the analytical solution of FODE \eqref{eq:fde} has the form of \eqref{eq:structure} when $f(t,u(t))$ satisfies some suitable conditions.
In the following, we will always assume that  $u(t)$ can be expressed in the following form for   convenience,
\begin{equation}\label{eq:analsol}
u(t)-u(0)=\sum_{r=1}^{m+1}c_rt^{\sigma_r}+\xi(t)t^{\sigma_{m+2}},
{\quad}0<\sigma_r<\sigma_{r+1},
\end{equation}
where $\xi(t)$ is a uniformly continuous function over the interval $[0,T]$ and $c_r\in\mathbb{R}$ are constants. { The sequence $\{\sigma_r\}$ is uniquely determined by the considered equation.  For example, when $f_u(t)=f(t)$ in \eqref{eq:fde} is smooth, $\sigma_r$'s are of the form $\{i+j\beta\}$, see
  \cite[Chapter 5]{Podlubny1999}. Another example is   from Lemma \ref{lem:diethm1},   $\sigma_r$'s are of the form \eqref{eq:structure}.}

Given a sequence of positive numbers $\{\theta_r\}$, we  define the operator
$I_{h,\theta}^{\beta,n,m}$  as
\begin{equation}\label{Ih}
I_{h,\theta}^{\beta,n,m}g=h^\beta \sum_{k=0}^n\omega_{n-k}^{(\beta)}g(t_k)
+h^\beta \sum_{k=1}^{m} W_{n,k}^{(\beta,\theta)} g(t_k)+h^\beta B_n^\theta g(t_0),
\end{equation}
where $\omega_{k}^{(\beta)}$ satisfies \eqref{eq:genfun}, and $W_{n,k}^{(\beta,\theta)}$ and $B_n^\theta$
are given by
\begin{eqnarray}
\sum_{k=1}^{m} W_{n,k}^{(\beta,\theta)}k^{\theta_r}&=&\frac{\Gamma(\theta_r+1)}{\Gamma(\theta_r+1+\beta)}n^{\theta_r+\beta}
-\sum_{k=0}^n\omega^{(\beta)}_{n-k}k^{\theta_r},{\quad}1\leq r \leq m,\label{eq:van1}\\
B_n^\theta&=&\frac{n^\beta}{\Gamma(1+\beta)}- \sum_{k=0}^n\omega_{n-k}^{(\beta)} -\sum_{k=1}^{m}{W}_{n,k}^{(\beta,\theta)}.\label{eq:bu}
\end{eqnarray}
Here $W_{n,k}^{(\beta,\theta)}$ in \eqref{Ih} are called the starting weights
that are chosen such that $I_{h,\theta}^{\beta,n,m}g=(I_0^\beta g)(t_n)$
when $g(t)=t^{\theta_r}\,(1\leq r \leq m)$, which leads to \eqref{eq:van1}.
%{
%\begin{rem}
%The linear system \eqref{eq:van1} is ill-conditioned when $m$ is
%large \cite{DieFFW06}. The large condition number of the Vandermonde-type matrix in \eqref{eq:van1} may lead to big roundoff errors of the starting weights $W_{n,k}^{(\beta,\theta)} (1\leq k\leq m)$.
%Fortunately, in practice we do not need many correction terms to get satisfactory numerical solutions,  so we can still have some reasonable accuracy for the starting weights. We present the residual of the system \eqref{eq:van1} in some numerical examples, which may determine, to some extent, the accuracy of \eqref{Ih} \cite{DieFFW06,Lub86}; see Table \ref{tbl:con1} for Example \ref{exm:linearsys1} and Table \ref{tbl:con2} for Example \ref{exm:nonlinear1} (Case II).  For more discussion on this issue, we refer to \cite{ZengZK2015}.
%\end{rem}}

{
\begin{rem}
The linear system \eqref{eq:van1} is ill-conditioned when $m$ is
large \cite{DieFFW06,Lub86,ZengZK2015}.  The large condition number of the Vandermonde-type matrix in \eqref{eq:van1} may lead to big roundoff errors of the starting weights $W_{n,k}^{(\beta,\theta)} (1\leq k\leq m)$.  However,
%the accuracy of the discrete operator $I_{h,\theta}^{\beta,n,m}$ depends on the number of the correction terms $m$ and the residue of the system \eqref{eq:van1} (see \cite{DieFFW06,Lub86}). Furthermore,
 we do not need many correction terms to get satisfactory numerical solutions in computation as observed in \cite{ZengZK2015}.  With this observation, we only need to solve the system \eqref{eq:van1}
with moderately large condition number. Thus we  can  obtain  reasonable accuracy of the staring weights
and hence the numerical solutions, see \cite{DieFFW06,ZengZK2015}.
We present   residuals of the system \eqref{eq:van1} and its corresponding condition numbers in    Example \ref{exm:linearsys1} and  Example \ref{exm:nonlinear1}.
%For more discussion on this issue, we refer to \cite{DieFFW06,ZengZK2015}.
\end{rem}}

If we apply \eqref{Ih} to approximate
$I_{0}^\beta  u(t)$, where $u(t)$ satisfies \eqref{eq:analsol},
then by \eqref{eq:lubich-1} we have
\begin{eqnarray}
(I_{0}^\beta u) (t_n)&=&(I_{0}^\beta (u-u(0)))(t_n)
+\frac{t_n^{\beta}u(0)}{\Gamma(1+\beta)}\nonumber\\
&=&h^\beta \sum_{k=0}^n\omega_{n-k}^{(\beta)} (u(t_k)-u_0)
+h^\beta \sum_{k=1}^{m_u} W_{n,k}^{(\beta,\sigma)} (u(t_k)-u_0)
+\frac{t_n^{\beta}u_0}{\Gamma(1+\beta)}+R_u^n\nonumber\\
&=&I_{h,\sigma}^{\beta,n,m_u}u+R_u^n,\label{eq:lubich}
\end{eqnarray}
where $I_{h,\sigma}^{\beta,n,m_u}$ is defined by \eqref{Ih} and $R^n_u$ is defined by
\begin{eqnarray}
R^n_u&=&O(h^2t_n^{\sigma_{m_u+1}+\beta-2})+O(h^{1+\sigma_{m_u+1}}t_n^{\beta-1}).
\label{error-ru}
\end{eqnarray}

From \eqref{eq:fde}, we have
\begin{eqnarray}\label{eq:nonlinear1}
f_u(t)=f(t,u(t))=(\,^CD_0^{\beta}u)(t)-\lambda u(t).
\end{eqnarray}
Hence, the regularity of $f(t,u(t))$ is  related to the regularity of $u$. In fact, based on the smoothness assumption of  $u(t)$ (see Eq. \eqref{eq:analsol}), we obtain that  $f(t,u(t))$ has the form
\begin{equation}\label{eq:analf}
\begin{aligned}
f(t,u(t))-f(0,u(0))=&-\lambda\sum_{r=1}^{m+1}c_rt^{\sigma_r}
+ \sum_{r=1}^{m+1}c_r\frac{\Gamma(\sigma_r+1)}{\Gamma(\sigma_r+1-\beta)}t^{\sigma_r-\beta}+\cdots\\
=&\sum_{r=1}^{l+1}d_rt^{\delta_r} + \zeta(t)t^{\delta_{l+2}},
\end{aligned}\end{equation}
where { $\zeta(t)$ is uniformly bounded on $[0,T]$, $\delta_r<\delta_{r+1}$, and
$\delta_r\in \{\sigma_k,k\geq 1\}\cup \{\sigma_k-\beta,k\geq 1\}$ with $\sigma_k$'s    from  \eqref{eq:analsol}.}

Similar to \eqref{eq:lubich}, we have
\begin{eqnarray}\label{approx-f}
(I_{t_0}^\beta f_u)(t_n)&=&I_{h,\delta}^{\beta,n,m_f}f_u +R_f^n,
\end{eqnarray}
where $I_{h,\delta}^{\beta,n,m_f}$ is defined by \eqref{Ih},

and the truncation error $R^n_f$ is given by
\begin{eqnarray}
%B_n^\delta &=&\frac{n^\beta}{\Gamma(1+\beta)}- \sum_{k=0}^n\omega_{n-k}^{(\beta)} -\sum_{k=1}^{m_f}{W}_{n,k}^{(\beta,f)},\label{eq:bf}\\
R_f^n&=&O(h^2t_n^{\delta_{m_f+1}+\beta-2})+O(h^{1+\delta_{m_f+1}}t_n^{\beta-1}). \label{eq:errf}
\end{eqnarray}

From \eqref{eq:lubich} and \eqref{approx-f}, we can derive the following implicit discretization for
\eqref{eq:int}
\begin{equation}\label{scheme}
u(t_n)=u_0+\lambda  I_{h,\sigma}^{\beta,n,m_u}u+I_{h,\delta}^{\beta,n,m_f}f_u+R_u^n+R_f^n,
\end{equation}
where $I_{h,\sigma}^{\beta,n,m_u}$ and  $I_{h,\delta}^{\beta,n,m_f}$ are defined by \eqref{Ih},
$R_u^n$ and $R_f^n$ are defined by  \eqref{error-ru} and \eqref{eq:errf}, respectively.

Let $U_k$ be the approximate solution of $u(t_k)$. Dropping the truncation errors $R^n_u$ and $R^n_f$ in
\eqref{scheme} and replacing $u(t_k)$ with $U_k$, we derive the fully implicit method for \eqref{eq:int}: to find $U_n$ for $n=n_0,n_0+1,...$ such that
\begin{equation}\label{im-scheme}
\begin{aligned}
U_n=&u_0+\lambda I_{h,\sigma}^{\beta,n,m_u}U+ I_{h,\delta}^{\beta,n,m_f}F,
\end{aligned}\end{equation}
where $n_0=1+\max\{m_f,m_u\}$, $F_n=f(t_n,U_n)$.

Given $U_k$  $(k=0,1,\ldots,n-1)$,  we need to solve a nonlinear system   \eqref{im-scheme}  at each time step to get $U_n$.
Next, we further use the extrapolation and Taylor expansion with correction terms to approximate $f(t_n,u(t_n))$ in \eqref{scheme}, which leads to linear systems and also preserves high-order accuracy.

%We need the following results:
If $u(t) = t^{\sigma}, \sigma>0$, then we have from the Taylor expansion that
  \begin{eqnarray}
  &&u(t_n)=2u(t_{n-1}) - u(t_{n-2})+O(h^2t_n^{\sigma-2}), \;\;{n\geq 2},\label{extrap-error}\\
  &&u(t_n)=u(t_{n-1})+hu'(t_{n-1}) +O(h^2t_n^{\sigma-2}),\;\;{n\geq 2}.\label{taylor-error}
\end{eqnarray}

%\begin{itemize}
(1) \textbf{By extrapolation with correction terms:}  It is clear that \eqref{extrap-error}
does not preserve globally second-order accuracy when $\sigma<2$.
Hence,   $2f(t_{n-1},u(t_{n-1}))-f(t_{n-2},u(t_{n-2}))$ is not a second-order approximation of $f(t_n,u(t_n))$ when $f(t,u(t))$ is not sufficiently smooth, see \eqref{eq:analf}.
By adding correction terms to  \eqref{extrap-error}, we can obtain
\begin{eqnarray}
f(t_n,u(t_n))&=&2f(t_{n-1},u(t_{n-1}))-f(t_{n-2},u(t_{n-2}))\nonumber\\
&&+\sum_{k=1}^{\widetilde{m}_f}\widehat{W}^{(f)}_{n,k}\left(f(t_k,u(t_k))-f(t_0,u_0)\right)
+\widetilde{R}_{f}^n,{\quad}n\geq 2,\label{eq:extra1}
\end{eqnarray}
where  $\{\widehat{W}^{(f)}_{n,k}\}$ are  chosen such that {the above equation \eqref{eq:extra1} is exact}, i.e., $\widetilde{R}_{f}^n=0$, for $f(t,u(t))=t^{\delta_r}(1\leq r\leq \widetilde{m}_f)$, i.e.,  $\{\widehat{W}^{(f)}_{n,k}\}$ satisfy
\begin{equation}\label{eq:extraweights}
\sum_{k=1}^{\widetilde{m}_f}\widehat{W}^{(f)}_{n,k} k^{\delta_r}=n^{\delta_r}-2({n-1})^{\delta_r}+({n-2})^{\delta_r},\quad r=1,\cdots,\widetilde{m}_f.
\end{equation}
The truncation error $\widetilde{R}_{f}^n$ in \eqref{eq:extra1} satisfies
\begin{equation}\label{eq:errc}
\widetilde{R}_{f}^n = O(h^2t_n^{\delta_{\widetilde{m}_f+1}-2})
\end{equation}
when  $f(t,u(t))$ satisfies \eqref{eq:analf}.

Inserting \eqref{eq:extra1}  into \eqref{scheme} yields
\begin{equation}\label{eq:intf1}\begin{aligned}
u(t_n)=&u_0+\lambda  I_{h,\sigma}^{\beta,n,m_u}u+I_{h,\delta}^{\beta,n,m_f}f_u\\
&+h^{\beta}\omega_0^{(\beta)}\bigg[-f(t_{n},u(t_{n}))+2f(t_{n-1},u(t_{n-1}))-f(t_{n-2},u(t_{n-2}))\\
&+\sum_{k=1}^{\widetilde{m}_f} \widehat{W}_{n,k}^{(f)} (f(t_k,u(t_k))-f(t_0,u_0))\bigg]+ R^n_E,
\end{aligned}
\end{equation}
where $R^n_E=R_u^n+R_f^n+h^{\beta}\widetilde{R}_{f}^n$, $R_u^n$, $R_f^n$, and $\widetilde{R}_{f}^n$
are defined by \eqref{error-ru}, \eqref{eq:errf}, and \eqref{eq:errc}, respectively.

From \eqref{eq:intf1}, we  obtain the IMEX method  based on the   extrapolation technique (abbreviated as \textbf{IMEX-E}) as:
given $U_k(0\leq k \leq n-1)$,  to find $U_n\;(n\geq 2)$ such that

\begin{equation}\label{scheme-ts-e}
\begin{aligned}
U_n=&U_0+\lambda I_{h,\sigma}^{\beta,n,m_u}U+ I_{h,\delta}^{\beta,n,m_f}F-h^{\beta}\omega_0^{(\beta)}F_n\\
&+h^{\beta}\omega_0^{(\beta)}\bigg[2F_{n-1}-F_{n-2}
+\sum_{k=1}^{\widetilde{m}_f} \widehat{W}_{n,k}^{(f)} (F_k-F_0)\bigg],
\end{aligned}\end{equation}
where $F_k=f(t_k,U_k)$,  $\omega_0^{(\beta)}=2^{-\beta}$, $I_{h,\sigma}^{\beta,n,m_u}$ and    $I_{h,\delta}^{\beta,n,m_f}$ are defined by \eqref{Ih}, and  $\widehat{W}_{n,k}^{(f)}$  is given by \eqref{eq:extraweights}.
\begin{rem}
Given $U_k(0\leq k\leq  n-1)$, Eq. \eqref{scheme-ts-e} is a linear equation of $U_n$. In fact,
$I_{h,\delta}^{\beta,n,m_f}F$ contains $h^{\beta}\omega_0^{(\beta)}F_n$, and it can be eliminated by the following term $-h^{\beta}\omega_0^{(\beta)}F_n$ in the scheme.
\end{rem}

(2) \textbf{By Taylor expansion with correction terms:} From the Taylor expansion and \eqref{taylor-error}, we   have
\begin{equation}\label{eq:zeng1}\begin{aligned}
f(t_n,u(t_n))=&f(t_{n-1},u(t_{n-1}))+hf'(t_{n-1},u(t_{n-1}))\\
&+\sum_{k=1}^{\widetilde{m}_f}\widetilde{W}^{(f)}_{n,k}\left(f(t_k,u(t_k))-f(t_0,u_0)\right)
+ O(h^2t_n^{\delta_{\widetilde{m}_f+1}-2}),\;{n\geq 2,}
\end{aligned}\end{equation}
where $f'(t_{n-1},u(t_{n-1}))=\frac{\mathrm{d}}{\mathrm{d}t}f(t,u(t))|_{t=t_{n-1}}$, and the starting weights $\{\widetilde{W}^{(f)}_{n,k}\}$ are chosen such that  Eq.  \eqref{eq:zeng1} is exact for $f(t,u(t)) = t^{\delta_r}\,
(1\leq r \leq  \widetilde{m}_f)$,
i.e., $\{\widetilde{W}^{(f)}_{n,k}\}$  satisfy
\begin{equation}\label{eq:zeng1-2}\begin{aligned}
\sum_{k=1}^{\widetilde{m}_f}\widetilde{W}^{(f)}_{n,k}k^{\delta_r}
=n^{\delta_r}-(n-1)^{\delta_r} - {\delta_r}(n-1)^{\delta_r-1},
{\quad} 1\leq r \leq  \widetilde{m}_f.
\end{aligned}\end{equation}
%\begin{equation}\label{eq:zeng1-2}\begin{aligned}
%u(t_n)=u(t_{n-1})+hu'(t_{n-1}))
%+\sum_{k=1}^{\widetilde{m}_f}\widetilde{W}^{(f)}_{n,k}u(t_k),
%\quad u(t) = t^{\delta_r}, r = 1,2,...,m_f.
%\end{aligned}\end{equation}

Next, we approximate $f'(t_{n-1},u(t_{n-1}))$ with correction terms, which is given by
\begin{eqnarray}\label{eq:zeng2}
f'(t_{n-1},u(t_{n-1}))&= &\px[u]f(t_{n-1},u(t_{n-1}))\px[t]u(t_{n-1})
+\px[t]f(t_{n-1},u(t_{n-1})) \notag \\
&=&\px[u]f(t_{n-1},u(t_{n-1}))\bigg[\frac{u(t_n)-u(t_{n-1})}{h}
+\frac{1}{h}\sum_{k=1}^{\widetilde{m}_{u}}\widetilde{W}^{(u)}_{n,k}(u(t_k)-u_0) \notag \\
&&+O(ht_n^{\sigma_{\widetilde{m}_u+1}-2})\bigg]+\px[t]f(t_{n-1},u(t_{n-1})),\;  n\geq 2,
 \end{eqnarray}
where $\{\widetilde{W}^{(u)}_{n,k}\}$ are chosen such that
$$\frac{u(t_n)-u(t_{n-1})}{h}
+\frac{1}{h}\sum_{k=1}^{\widetilde{m}_{u}}\widetilde{W}^{(u)}_{n,k}(u(t_k)-u_0)=u'(t_{n-1})$$
for some $u(t)=t^{\sigma_r}(1\leq r \leq \widetilde{m}_u)$, i.e., $\{\widetilde{W}^{(u)}_{n,k}\}$
satisfy
\begin{equation}\label{eq:zeng2-2}\begin{aligned}
\sum_{k=1}^{\widetilde{m}_{u}}\widetilde{W}^{(u)}_{n,k}k^{\sigma_r}
=\sigma_r(n-1)^{\sigma_r-1}- (n^{\sigma_r} - (n-1)^{\sigma_r}),
{\quad}1\leq r \leq \widetilde{m}_u.
\end{aligned}\end{equation}
Combining \eqref{eq:zeng1} and \eqref{eq:zeng2} yields
\begin{eqnarray}\label{eq:taylor1}
f(t_n,u(t_n))&=&f(t_{n-1},u(t_{n-1}))
+h\px[t]f(t_{n-1},u(t_{n-1})) \notag \\
&&+\px[u]f(t_{n-1},u(t_{n-1}))\left[u(t_n)-u(t_{n-1})
+\sum_{k=1}^{\widetilde{m}_u}\widetilde{W}^{(u)}_{n,k}(u(t_k)-u_0)\right] \notag\\
&&+\sum_{k=1}^{\widetilde{m}_f}\widetilde{W}^{(f)}_{n,k}\left(f(t_k,u(t_k))-f(t_0,u_0)\right) +\widetilde{R}_{u,f}^n,
 \end{eqnarray}
where the truncation error $\widetilde{R}_{u,f}^n$ satisfies
\begin{equation}\label{eq:taylorerror}
\begin{aligned}
\widetilde{R}_{u,f}^n= &O(h^2t_n^{\delta_{\widetilde{m}_f+1}-2})
+O(h^2t_n^{\sigma_{\widetilde{m}_u+1}-2}).
\end{aligned}\end{equation}
Inserting \eqref{eq:taylor1}  into \eqref{scheme} leads to
\begin{equation}\label{eq:intf2}
\begin{aligned}
u(t_n)=&u_0+\lambda I_{h,\sigma}^{\beta,n,m_u}u+ I_{h,\delta}^{\beta,n,m_f}f_u\\
&+h^{\beta}\omega_0^{(\beta)}\bigg[-f(t_{n},u(t_{n}))+f(t_{n-1},u(t_{n-1}))
+h\px[t]f(t_{n-1},u(t_{n-1}))\\
&+\px[u]f(t_{n-1},u(t_{n-1}))\Big(u(t_n)-u(t_{n-1})
+\sum_{k=1}^{\widetilde{m}_u}\widetilde{W}^{(u)}_{n,k}(u(t_k)-u_0)\Big)\\
&+\sum_{k=1}^{\widetilde{m}_f}\widetilde{W}^{(f)}_{n,k}
\left(f(t_k,u(t_k))-f(t_0,u_0)\right)\bigg]+ R^n_T,
\end{aligned}
\end{equation}
where $R^n_T=R_u^n+R_f^n+h^{\beta}\widetilde{R}_{u,f}^n$,
$R_u^n$, $R_f^n$, and $\widetilde{R}_{u,f}^n$
are defined by \eqref{error-ru}, \eqref{eq:errf}, and \eqref{eq:taylorerror}, respectively.

From \eqref{eq:intf2}, we  obtain the IMEX method  based on the  Taylor expansion technique   (abbreviated as \textbf{IMEX-T}) as:
given $U_k(0\leq k \leq n-1)$,   to find $U_n\,(n\geq 2)$ such that
\begin{equation}\label{scheme-ts-t}
\begin{aligned}
U_n&=U_0+\lambda I_{h,\sigma}^{\beta,n,m_u}U+ I_{h,\delta}^{\beta,n,m_f}F
+h^{\beta}\omega_0^{(\beta)}\bigg[-F_n+F_{n-1}
+h\px[t]f(t_{n-1},U_{n-1})\\
&\;\;\,+\px[u]f(t_{n-1},U_{n-1})\Big(U_n-U_{n-1}
+\sum_{k=1}^{\widetilde{m}_u}\widetilde{W}^{(u)}_{n,k}(U_k-U_0)\Big)+\sum_{k=1}^{\widetilde{m}_f}\widetilde{W}^{(f)}_{n,k}\left(F_k-F_0\right)\bigg],
\end{aligned}\end{equation}
where $F_k=f(t_k,U_k)$,  $\omega_0^{(\beta)}=2^{-\beta}$, $I_{h,\sigma}^{\beta,n,m_u}$ and    $I_{h,\delta}^{\beta,n,m_f}$ are defined by \eqref{Ih},
$\widetilde{W}^{(f)}_{n,k}$ and  $\widetilde{W}^{(u)}_{n,k}$ are  given by
\eqref{eq:zeng1-2} and \eqref{eq:zeng2-2}, respectively.

%\begin{equation}\label{scheme-ts-t}
%\begin{aligned}
%U_n&=U_0+\lambda h^{\beta}\left[\sum_{k=0}^n\omega_{n-k}^{(\beta)}U_k
%+\sum_{k=1}^{m_u} W_{n,k}^{(\beta,\sigma)} U_k+B_n^\sigma U_0\right]\\
%&\;\;\,+h^{\beta}\left[\sum_{k=0}^{n-1}\omega_{n-k}^{(\beta)}F_k
%+\sum_{k=1}^{m_f} W_{n,k}^{(\beta,\delta)} F_k+B_n^\delta F_0\right]+h^{\beta}\omega_0^{(\beta)}\bigg[F_{n-1}
%+h\px[t]f(t_{n-1},U_{n-1})\\
%&\;\;\,+\px[u]f(t_{n-1},U_{n-1})\Big(U_n-U_{n-1}
%+\sum_{k=1}^{\widetilde{m}_u}\widetilde{W}^{(u)}_{n,k}(U_k-U_0)\Big)+\sum_{k=1}^{\widetilde{m}_f}\widetilde{W}^{(f)}_{n,k}\left(F_k-F_0\right)\bigg],
%\end{aligned}\end{equation}
%where $F_k=f(t_k,U_k)$,  $\omega_0^{(\beta)}=2^{-\beta}$,  $W_{n,k}^{(\beta,\sigma)}$ and   $W_{n,k}^{(\beta,\delta)}$ are defined by \eqref{eq:van1},
%$\widetilde{W}^{(f)}_{n,k}$ and  $\widetilde{W}^{(u)}_{n,k}$ are  given by
%\eqref{eq:zeng1-2} and \eqref{eq:zeng2-2}, respectively.

Next, we present the convergence results for the two schemes \eqref{scheme-ts-e} and \eqref{scheme-ts-t}, the proofs of which will be given in Section \ref{sec:proof}.
\begin{thm}[Convergence of  IMEX-E]\label{thm:convergence}
Suppose that  $u(t)$ is the  solution to \eqref{eq:fde} that satisfies \eqref{eq:analsol} and $U_n$ is the solution to \eqref{scheme-ts-e}, and $f(t,u)$ satisfies the Lipschitz condition with respect to the second argument $u$. If $\sigma_{m_u},\delta_{m_f}\leq 2$ and $\delta_{\widetilde{m}_f}\leq 2+\beta$,
then there exists a positive constant $C$ independent of $h$ and $n$ such that
\begin{equation}\label{eq:convergence}
|u(t_n)-U_n|\leq C\left(\sum_{k=1}^{m}|u(t_k)-U_k|
+h^{q}\right),
\end{equation}
where  $m=\max\{m_u,m_f,\widetilde{m}_f\}$ and
$
q=\min\{2,\sigma_{m_u+1}+\beta,\delta_{m_f+1}+\beta,\delta_{\widetilde{m}_f+1}+\beta\}.
$ %\end{eqnarray*}

\end{thm}

\begin{thm}[Convergence of  IMEX-T]\label{thm:convergence2}
Suppose that  $u(t)$ is the  solution to \eqref{eq:fde} that satisfies \eqref{eq:analsol} and $U_n$ is the solution to \eqref{scheme-ts-t}, and $f(t,u)$ satisfies the Lipschitz condition with respect to the second argument $u$. If $\sigma_{m_u},\delta_{m_f}\leq 2$ and $\sigma_{\widetilde{m}_u},\delta_{\widetilde{m}_f}\leq 2+\beta$, then there exists a positive constant $C$ independent of $h$ and $n$ such that
\begin{equation}\label{eq:convergence2}
|u(t_n)-U_n|\leq C\left(\sum_{k=1}^{m}|u(t_k)-U_k|+h^q\right),
\end{equation}
where $m=\max\{m_u,m_f,\widetilde{m}_f,\widetilde{m}_u\}$ and
\begin{eqnarray*}
q=\min\{2,\sigma_{m_u+1}+\beta,\delta_{m_f+1}+\beta,\delta_{\widetilde{m}_f+1}+\beta,
\sigma_{\widetilde{m}_u+1}+\beta\}.
\end{eqnarray*}
	
\end{thm}

%Note that the last term in the right side of \eqref{eq:convergence2} comes from \eqref{eq:taylorerror}, which will cause low accuracy at the origin when the exact solution is not smooth.  The other parts of the proof is very similar to the proof of Theorem \ref{thm:convergence}, so we omit it in this paper.
\section{Linear stability of IMEX schemes}\label{sec:stability-linear}

In this section, we discuss the  linear stability of the proposed IMEX schemes for the scalar equation
\begin{equation}\label{eq:testeq1}
(\,^CD_0^{\beta}u)(t)=\lambda u(t)+\rho u(t),\; t\geq 0,\; \lambda,\,\rho\in \mathbb{C}.
\end{equation}

We recall the definition of stability  for the linear equation  \eqref{eq:testeq1}.
\begin{thm}[\cite{Lub85,Lub86}]
	Let $\beta>0$. The steady-state solution $u=0$ of   Eq. {\eqref{eq:testeq1}}  is stable if and only if $(\lambda+\rho)\in\sum_\beta$, where $\sum_\beta=\{s\in\mathbb{C}:\;|arg(s)|>\frac{\beta\pi}{2}\}$.
\end{thm}

\begin{defn}
	A numerical method is said to be $A(\frac{\beta\pi}{2})$-stable if    its stability region for \eqref{eq:testeq1} contains the whole sector $\sum_\beta$.
\end{defn}

The following theorem is useful to determine stability regions  of the numerical schemes.

\begin{thm}[\cite{Gar10,Garrappa2015,Lub86-1}]\label{stable region}
	Let $\beta>0$.  Assume  that the sequence $\{g_n\}$ is convergent and that the quadrature weights $w_n$ ($n\geq 1$) satisfy
	\begin{eqnarray}
	w_n=\frac{n^{\beta-1}}{\Gamma(\beta+1)}+v_n,\;\;\sum_{n=1}^{\infty}|v_n|<\infty,\label{eq:series}
	\end{eqnarray}
	then the stability region of the convolution quadrature
	$
	y_n=g_n+\xi\sum_{j=0}^n w_{n-j}y_j\label{eq:series-2}
	$
	is
\begin{equation*}
\Sigma_\beta^{{Num}}=\big\{ \xi \in\mathbb{C}\big|1-\xi w^\beta(z)\neq 0:\;|z|\leq 1\big\},\;w^\beta(z)=\sum_{n=0}^\infty w_n z^n,
\end{equation*}
where $\xi=\lambda h^\beta$ or $\xi$ is some function of $\lambda h^\beta$.
\end{thm}

We first consider the linear stability of the IMEX-E  scheme for the test equation  \eqref{eq:testeq1}.
From the IMEX-E scheme \eqref{scheme-ts-e}, we  get
\begin{eqnarray}\label{eq:linearstab1}
U_n &=& U_0+{(\lambda+\rho)}h^{\beta}\left[\sum_{k=0}^n\omega_{n-k}^{(\beta)}U_k
+\sum_{k=1}^{m_u} W_{n,k}^{(\beta,\sigma)} U_k+B_n^\sigma U_0\right] \\
&&+{\rho}h^{\beta}\omega_0^{(\beta)}\left[-U_n+2U_{n-1}-U_{n-2}
+\sum_{k=1}^{\widetilde{m}_f} \widehat{W}_{n,k}^{(f)} (U_k-U_0)\right] \notag \\
&=&U_0+{(\lambda+\rho)}h^{\beta}\sum_{k=0}^n\omega_{n-k}^{(\beta)}U_k
-{\rho}h^{\beta}\omega_0^{(\beta)}\left(U_n-2U_{n-1}+U_{n-2}\right)
+h^{\beta}\sum_{k=0}^{m} {\omega}_{n,k}U_k,\notag
 \end{eqnarray}
where $\sum_{k=0}^{m} {\omega}_{n,k}U_k=(\lambda+\rho)\sum_{k=1}^{m_u} W_{n,k}^{(\beta,\sigma)} U_k+B_n^\sigma U_0+\rho \omega_0^{(\beta)}\sum_{k=1}^{\widetilde{m}_f} \widehat{W}_{n,k}^{(f)} (U_k-U_0)$. By comparing coefficients on both sides of the above identity, we can get ${\omega}_{n,k}$. Here we do not give the exact expression of ${\omega}_{n,k}$, since it does not
affect the stability analysis.

Denote $U(z)=\sum_{k=0}^{\infty}U_kz^k,|z|\leq 1$, $m_0=\max\{m_u,\widetilde{m}_f\}$. Then we have from \eqref{eq:linearstab1} that
\begin{equation*}\label{eq:linearstab1-2}
\begin{aligned}
\sum_{n=2}^{\infty}U_nz^n
=&U_0\sum_{n=2}^{\infty}z^n+{(\lambda+\rho)}h^{\beta}
\sum_{n=2}^{\infty}\left(\sum_{k=0}^n\omega_{n-k}^{(\beta)}U_k\right)z^n\\
&-{\rho}h^{\beta}\omega_0^{(\beta)}\sum_{n=2}^{\infty}\left(U_n-2U_{n-1}+U_{n-2}\right)z^n
+h^\beta\sum_{k=0}^{m_0}U_k \left(\sum_{n=2}^{\infty}{\omega}_{n,k}z^n\right),
\end{aligned}\end{equation*}
which leads to
\begin{equation}\label{eq:linearstab1-3}
\begin{aligned}
U(z)&-U_0-U_1z=U_0(1-z)^{-1}-U_0(1+z)\\
&+{(\lambda+\rho)}h^{\beta} \left(U(z)\omega^{(\beta)}(z)-\omega_{0}^{(\beta)}U_0
-(\omega_{0}^{(\beta)}U_1+\omega_{1}^{(\beta)}U_0)z\right)\\
&-{\rho}h^{\beta}\omega_0^{(\beta)} \left((1-z)^2U(z)-(1-2z)U_0-U_1z\right)
+h^\beta\sum_{k=0}^{m_0} {\tilde{\omega}}_{k}(z)U_k,
\end{aligned}\end{equation}
where ${\tilde{\omega}}_{k}(z)=\sum_{n=2}^{\infty}{\omega}_{n,k}z^n$. We simplify \eqref{eq:linearstab1-3} as
\begin{equation}\label{eq:linearstab1-4}
\begin{aligned}
\left(1-{(\lambda+\rho)}h^{\beta}\omega^{(\beta)}(z)+{\rho}h^{\beta}\omega_0^{(\beta)}(1-z)^2\right)U(z)
=&H(z)=\sum_{n=0}^{\infty}H_nz^n,
\end{aligned}\end{equation}
where $H(z)=U_0(1-z)^{-1}+(U_1-U_0)z-{(\lambda+\rho)}h^{\beta}
\left(\omega_{0}^{(\beta)}U_0+(\omega_{0}^{(\beta)}U_1+\omega_{1}^{(\beta)}U_0)z\right)
+{\rho}h^{\beta}\omega_0^{(\beta)} \left((1-2z)U_0+U_1z\right)
+h^{\beta}\sum_{k=0}^{m_0} {\tilde{\omega}}_{k}(z)U_k$.

It is readily  verified that $\{H_n\}$
is a convergent sequence if $\sigma_{m_u},\delta_{\widetilde{m}_f}\leq 2-\beta$; see Lemmas \ref{lem:weights} and \ref{lem:weight2}.
From Theorem \ref{stable region}, we obtain that the scheme \eqref{scheme-ts-e} is stable  if
\begin{equation}
1-(\lambda+\rho)h^\beta\omega^{(\beta)}(z)+\rho h^\beta\omega^{(\beta)}_0(1-z)^2\neq 0,
\quad\forall |z|\leq 1.
\end{equation}
We summarize the above argument and   have the following theorem.
\begin{figure}[!t]
	\caption{ Stability region of the IMEX-E scheme (shaded) for the test equation \eqref{eq:testeq1}; left: stability region varying with different $\beta$ and  $\rho=0.5\lambda$; right: stability region varying with different $\rho$ restricted to $\lambda+\rho=\mu$, and $\beta=0.2$.}
	\label{fig:ts_e1}
		\centering
	\epsfig{figure=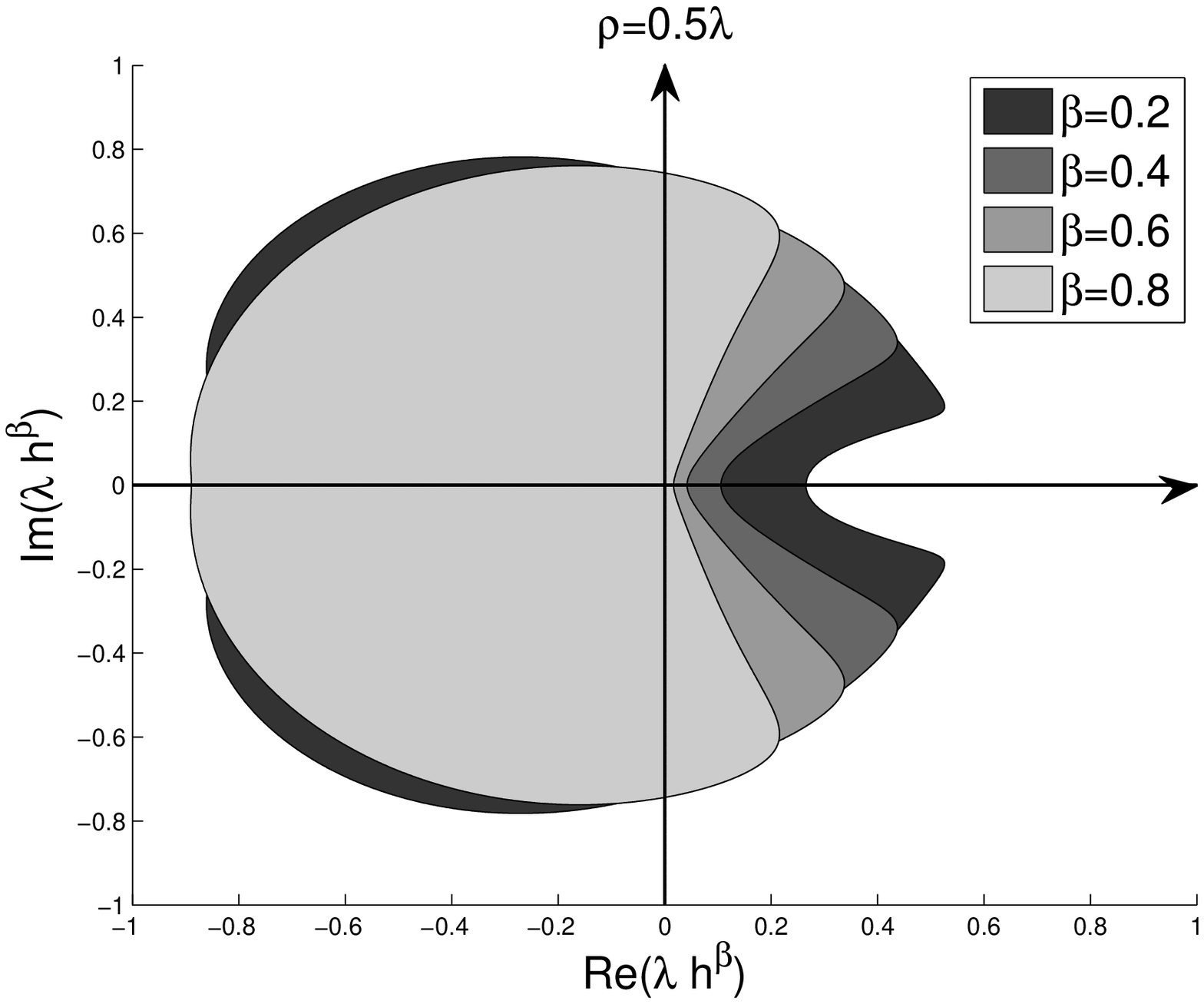,width=6cm}~\epsfig{figure=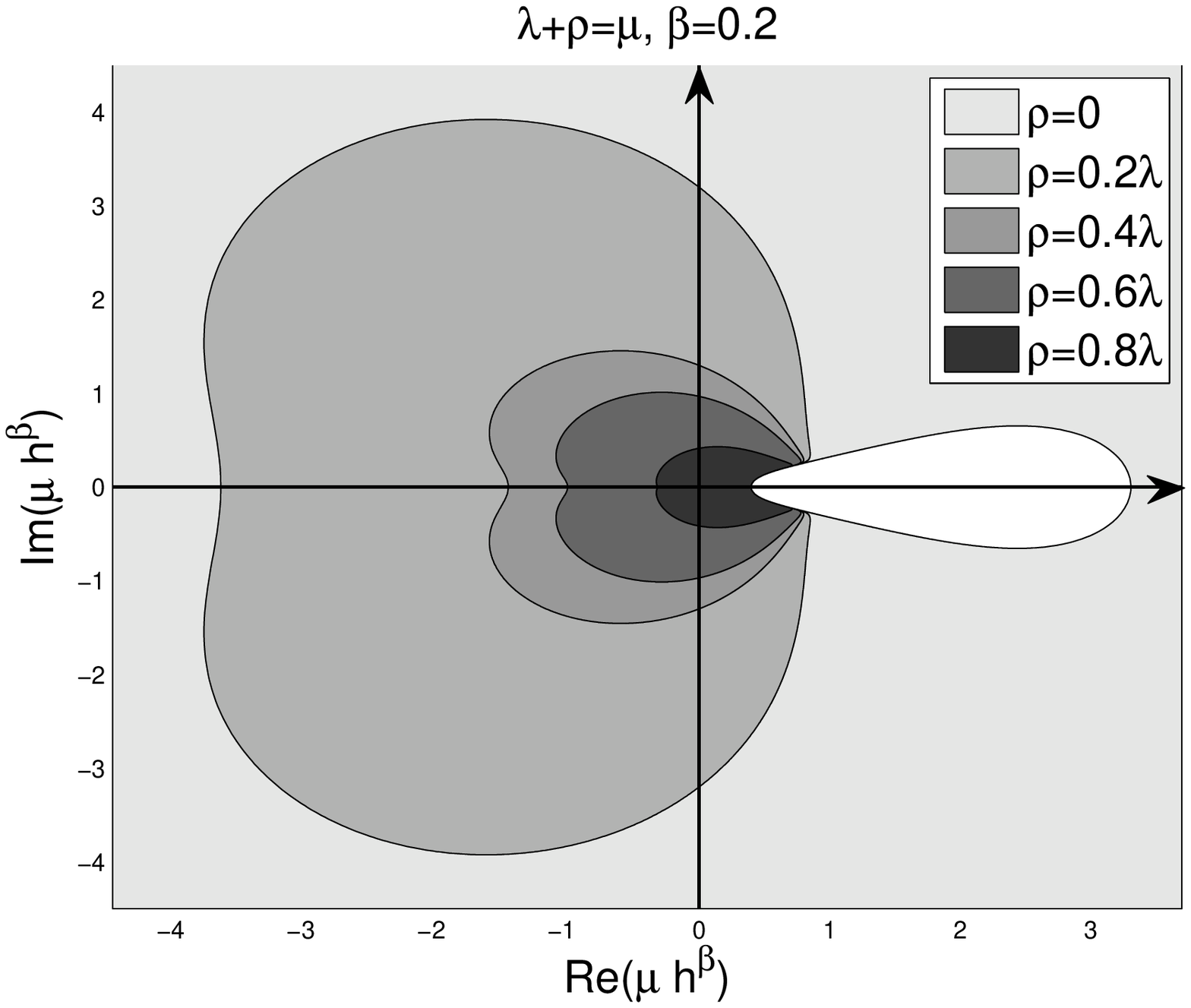,width=6cm}
\end{figure}
\begin{thm}[Linear stability of IMEX-E]\label{thm:stability1}
Let $\rho=k\lambda$ and
 $\xi=\lambda h^\beta$. Then the stability region of
the method \eqref{scheme-ts-e} for \eqref{eq:testeq1}  is
\begin{equation}
\Sigma_\beta^{Num}=\mathbb{C}\setminus \left\{\xi\Big|\xi=\frac{1}{(k+1)\omega^{(\beta)}(z)-k \omega^{(\beta)}_0(1-z)^2},\,|z|\leq 1 \right\},
\end{equation}	
\end{thm}

Fig. \ref{fig:ts_e1} shows the   stability region of the IMEX-E scheme with different $\beta$  and $\rho/\lambda$. Moreover, if $\lambda,\rho\in \mathbb{R}$ and  $\rho+\lambda=constant$, i.e., $\rho+\lambda=-1$, then we have $-1\leq \rho \leq 0$. In such a case, we should have  $h^{\beta}\neq \frac{2^{\beta}(1-z)^{\beta}}{-(1+z)^{\beta}-\rho (1-z)^{2+\beta}},|z|\leq 1$. Hence, the stability condition of the method \eqref{scheme-ts-e} for the model problem \eqref{eq:testeq1} with $\lambda+\rho=-1$ satisfies
$$h^{\beta}\in \left(0,-\frac{2^{\beta}}{4\rho}\right).$$
Clearly,  the length of the stability interval  decreases  as $\rho\to-1$, which is consistent with the theoretical result of the case $\beta=1$.

\begin{figure}[!t]
	\caption{ Stability region of the IMEX-T scheme (shaded) for the test equation \eqref{eq:testeq1}; left: $\beta=0.2$, right: $\beta=0.8$}
	\label{fig:ts_t1}
	\centering
\epsfig{figure=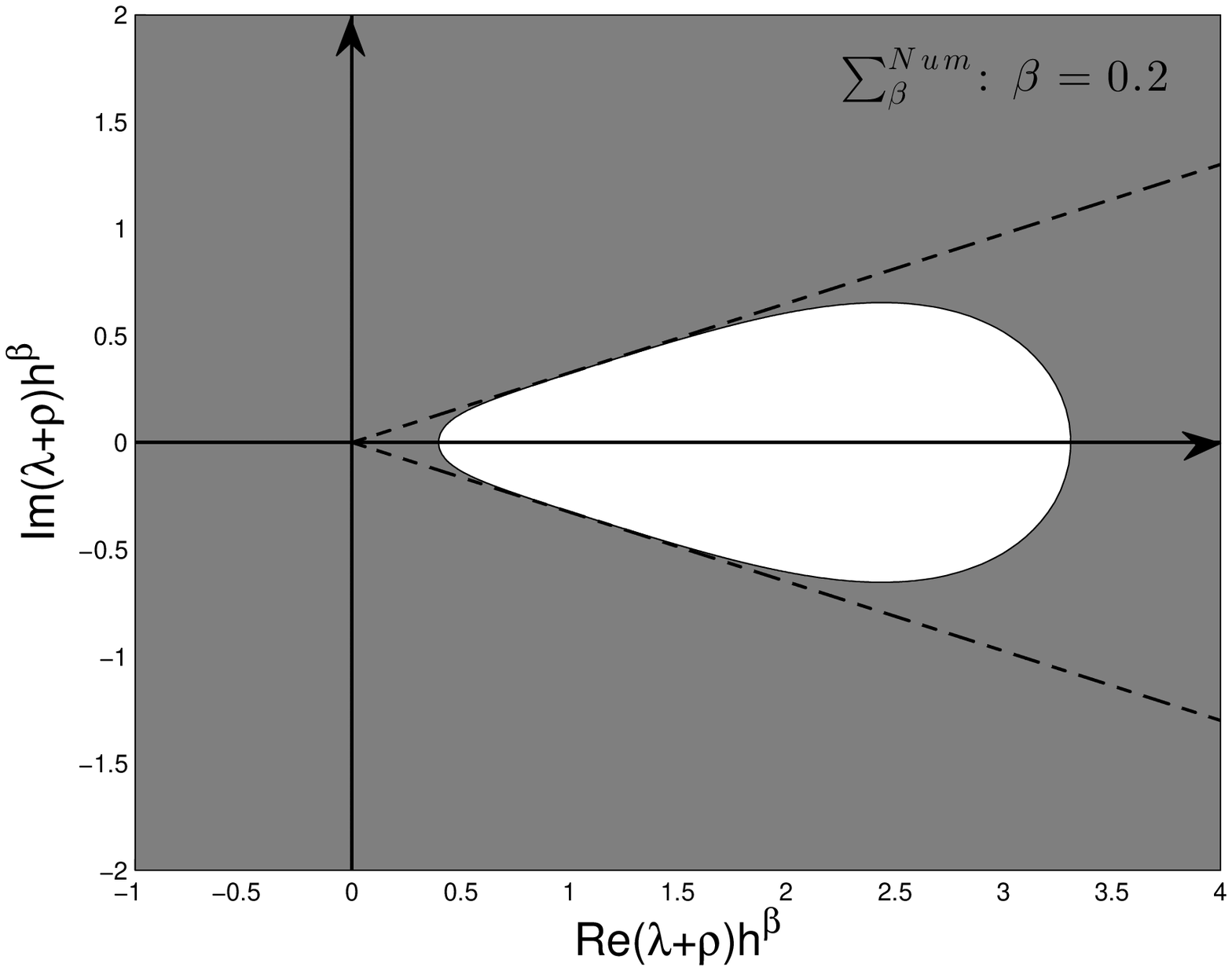,width=6cm}~\epsfig{figure=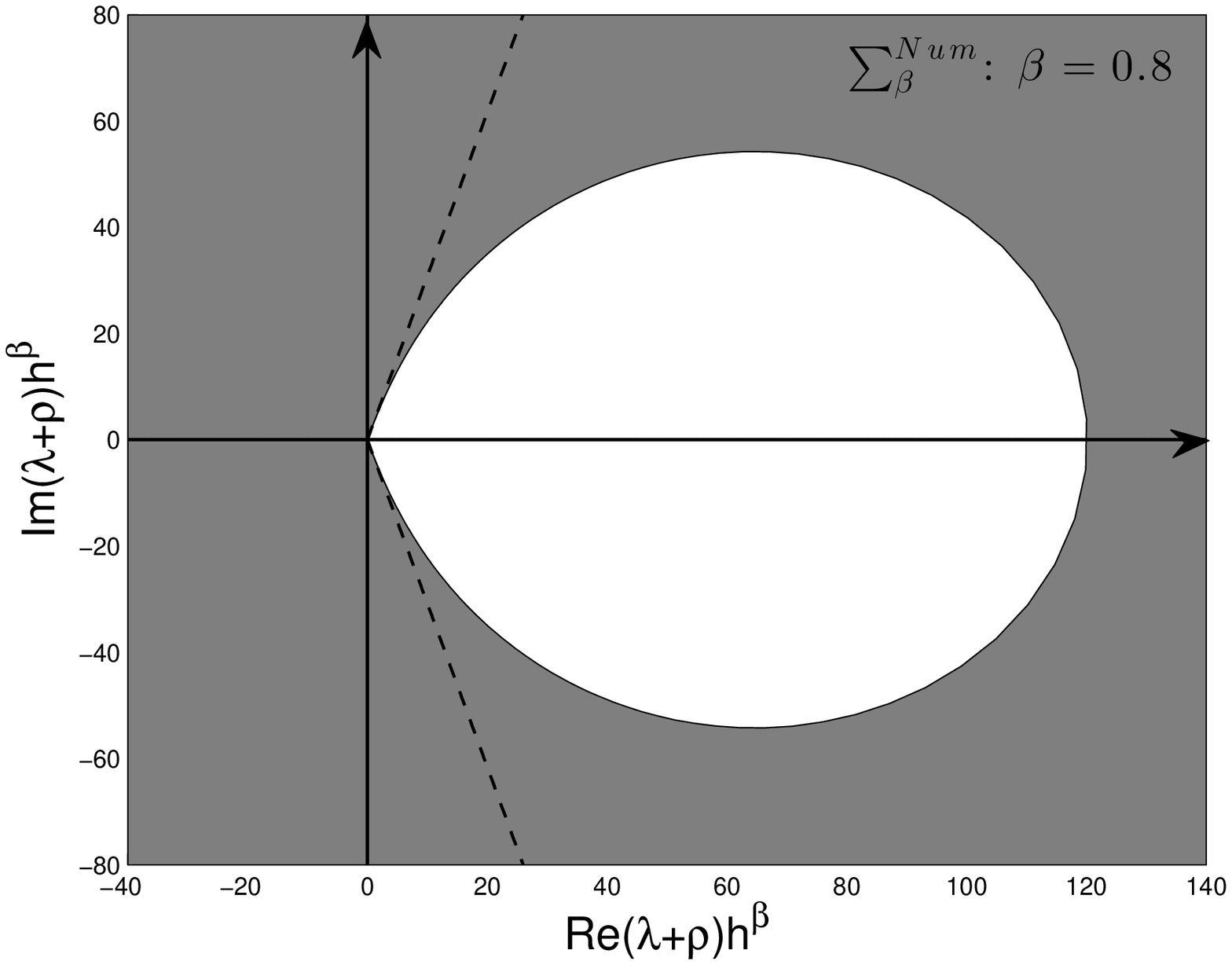,width=6cm}
\end{figure}

Next we consider the linear stability of the scheme IMEX-T \eqref{scheme-ts-t}.  Applying \eqref{scheme-ts-t} to \eqref{eq:testeq1} and letting $f(t,U) = (\rho+\lambda) U$, we have
the stability region of \eqref{scheme-ts-t} in the following theorem, see also \cite{Lub86-1} and Fig. \ref{fig:ts_t1}.

\begin{thm}[Linear stability of IMEX-T]\label{thm:stability2}
The stability region of
the method \eqref{scheme-ts-t} for \eqref{eq:testeq1}  is
\begin{equation}
\Sigma_\beta^{Num}=\mathbb{C}\setminus \left\{\xi\bigg|
\xi=\frac{1}{\omega^{(\beta)}(z)}, |z|\leq 1\right\},
\end{equation}	
where $\xi=(\lambda+\rho) h^\beta$, and $\omega^{(\beta)}(z)$ is  defined by \eqref{eq:genfun}.
\end{thm}

\section{Extensions}\label{sec4}
Besides the IMEX-E and IMEX-T schemes that we have derived, some other second-order IMEX schemes  can be obtained by using different approximations of $(I_0^\beta u)(t_n)$.

Let $\{\omega_{j}^{(\beta)}\}$ be coefficients of the Taylor expansion of the following  generating function
\begin{eqnarray}
\omega^{(\beta)}(z)&=&(1-z)^{-\beta}\left(1-\frac{\beta}{2}(1-z)\right)=\sum_{j=0}^\infty
\omega_{j}^{(\beta)}z^j, \label{eq:genfun11}\\
\text{or}\qquad\qquad\quad\qquad\omega^{(\beta)}(z)&=&(3/2-2z-z^2/2)^{-\beta}=\sum_{j=0}^\infty\omega_{j}^{(\beta)}z^j.\label{eq:genfun2}\qquad\qquad\qquad\qquad\qquad
\end{eqnarray}
Replacing the generating function   \eqref{eq:genfun}   with \eqref{eq:genfun11} or \eqref{eq:genfun2}  and  repeating the procedures below Eq. \eqref{eq:genfun} in Section \ref{sec:time-splitting-schemes} leads to new IMEX schemes, which have exactly the same form as the IMEX-E and IMEX-T schemes but using different weights.

Next, we present the trapezoidal rule \cite{DieFF04}  with correction terms that is given by
\begin{eqnarray}
I_{0}^{\beta}u(t_n)&=&\left[I_{0}^{\beta}(u(t)-u(t_0))\right]_{t=t_n}
+\frac{u(t_0)t_n^{\beta}}{\Gamma(1+\beta)}\nonumber\\
&=&h^\beta\sum_{j=0}^n b_{n,j}^{(\beta)} (u(t_j)-u_0) +{h^\beta} \sum_{j=1}^{m_u}\widetilde{W}_{n,j}^{(\beta,\sigma)}(u(t_j) -u_0) +\frac{u_0t_n^{\beta}}{\Gamma(1+\beta)}+O(h^2)\nonumber\\
&=&h^\beta\sum_{j=1}^n b_{n,j}^{(\beta)} u(t_j)
+{h^\beta} \sum_{j=1}^{m_u}\widetilde{W}_{n,j}^{(\beta,\sigma)}u(t_j)
+h^{\beta}\widetilde{B}_n^\sigma u_0+ O(h^2),\label{eq:trapapp}
\end{eqnarray}
where $b_{n,0}^{(\beta)}=[(n-1)^{\beta+1}-(n-1-\beta)n^{\beta}]/\Gamma(2+\beta)$, $b_{n,n}^{(\beta)}=1/\Gamma(2+\beta)$,
$b_{n,j}^{(\beta)}=\frac{1}{\Gamma(2+\beta)}
[(n-j+1)^{\beta+1}-2(n-j)^{\beta+1}+(n-j-1)^{\beta+1}]\,(1\leq j \leq n-1),$
and the starting weights  $\{\widetilde{W}_{n,j}^{(\beta,\sigma)} \}$ can be derived by solving the following linear system
\begin{equation}\label{eq:van2}
\sum_{k=1}^{m_u} \widetilde{W}_{n,k}^{(\beta,\sigma)}k^{\sigma_r}=\frac{\Gamma(\sigma_r+1)}{\Gamma(\sigma_r+1+\beta)}n^{\sigma_r+\beta}
-\sum_{k=1}^nb^{(\beta)}_{n,k}k^{\sigma_r},{\quad}1\leq r \leq m_u,
\end{equation}
and $\widetilde{B}_n^\sigma$ is given by
\begin{equation*}
	\widetilde{B}_n^\sigma=\frac{n^\beta}{\Gamma(1+\beta)}- \sum_{k=1}^nb_{n,k}^{(\beta)} -\sum_{k=1}^{m_u}\widetilde{W}_{n,k}^{(\beta,\sigma)}.
	\end{equation*}

Similar to deriving the IMEX-E scheme \eqref{scheme-ts-t}, we get the IMEX scheme for FODE \eqref{eq:fde} based on the extrapolation and the trapezoidal rule   (abbreviated as \textbf{IMEX-E-Trap}): given $U_k(0\leq k \leq n-1)$,  to find $U_n\;(n\geq 2)$ such that
\begin{equation}\label{scheme-ts-trap}
\begin{aligned}
U_n=&U_0+\lambda h^{\beta}\left[\sum_{k=1}^n b_{n,k}^{(\beta)}U_k
+\sum_{k=1}^{m_u} \widetilde{W}_{n,k}^{(\beta,\sigma)} U_k+\widetilde{B}_n^\sigma U_0\right]\\
&+h^{\beta}\left[\sum_{k=1}^{n-1}b_{n,k}^{(\beta)}F_k
+\sum_{k=1}^{m_f} \widetilde{W}_{n,k}^{(\beta,\delta)} F_k+\widetilde{B}_n^\delta F_0\right]\\
&+h^{\beta}b_{n,n}^{(\beta)}\bigg[2F_{n-1}-F_{n-2}
+\sum_{k=1}^{\widetilde{m}_f} \widehat{W}_{n,k}^{(f)} (F_k-F_0)\bigg],
\end{aligned}\end{equation}
where   $F_n=f(t_n,U_n)$, $b_{n,k}^{(\beta)}$ is given in \eqref{eq:trapapp}, $\widetilde{W}_{n,k}^{(\beta,\sigma)}$ is defined in \eqref{eq:van2},  $\widehat{W}^{(f)}_{n,k}$ is given by \eqref{eq:extraweights}, $\widetilde{B}_n^\delta=\frac{n^\beta}{\Gamma(1+\beta)}- \sum_{k=1}^nb_{n,k}^{(\beta)} -\sum_{k=1}^{m_f}\widetilde{W}_{n,k}^{(\beta,\delta)}$, and
$\widetilde{W}_{n,k}^{(\beta,\delta)}$ can be derived by solving the linear system
\begin{equation}\label{eq:van3}
\sum_{k=1}^{m_f} \widetilde{W}_{n,k}^{(\beta,
	\delta)}k^{\delta_r}=\frac{\Gamma(\delta_r+1)}{\Gamma(\delta_r+1+\beta)}n^{\delta_r+\beta}
-\sum_{k=1}^nb^{(\beta)}_{n,k}k^{\delta_r},{\quad}1\leq r \leq m_f.
\end{equation}

The Taylor expansion used in the IMEX-T scheme can be also applied here to get another IMEX scheme with the trapezoidal rule. We omit the details here due to the similarity.

%\subsection{Connection to splitting schemes}
The present IMEX-E and IMEX-T schemes can be extended to
the following multi-term fractional ordinary differential system
\begin{equation}\label{eq:multifde}
(\,^CD_{0}^{\alpha}u)(t)+(\,^CD_{0}^{\beta}u)(t)=Au(t)+f_u(t),\; t\in(0,T],\; u(0)=u_0,
 \end{equation}
where $A$ is a real-valued matrix. Let $\alpha<\beta$. We can transform  \eqref{eq:multifde} into its integral form as
\begin{eqnarray}\label{eq:mulint}
{u}(t)&=&{u}(0)-(I_{0}^{\beta-\alpha} (u-u(0))(t)+A\,(I_{0}^\beta u)(t) +(I_{0}^\beta f_u)(t),\quad 0\leq t\leq T.
\end{eqnarray}
Then we apply \eqref{Ih} to discretize each fractional integral in \eqref{eq:mulint}, and corresponding  IMEX schemes similar to \eqref{scheme-ts-e}, \eqref{scheme-ts-t} and \eqref{scheme-ts-trap} can be derived. For instance,
the IMEX-E scheme for \eqref{eq:multifde} reads
\begin{equation}\label{scheme-ts-e-mul}
\begin{aligned}
U_n=&U_0-I_{h,\sigma}^{\beta-\alpha,n,m_u^1}(U-U_0)
+AI_{h,\sigma}^{\beta,n,m_u^2}U+I_{h,\delta}^{\beta,n,m_f}F\\
&+h^{\beta}\omega_0^{(\beta)}\bigg[-F_n+2F_{n-1}-F_{n-2}
+\sum_{k=1}^{\widetilde{m}_f} \widehat{W}_{n,k}^{(f)} (F_k-F_0)\bigg],
\end{aligned}\end{equation}
where $I_{h,\sigma}^{\beta-\alpha,n,m_u^1}$, $I_{h,\sigma}^{\beta,n,m_u^2}$ and $I_{h,\delta}^{\beta,n,m_f}$  are defined by \eqref{Ih} and $\widehat{W}_{n,k}^{(f)}$ is defined by \eqref{eq:extraweights}.
Since the IMEX-T and IMEX-Trap schemes can be also derived readily, we do not present them here.
%\begin{equation}\label{scheme-ts-e-mul}
%\begin{aligned}
%U_n=&U_0-h^{\beta-\alpha}\left[\sum_{k=0}^n\omega_{n-k}^{(\beta-\alpha)}U_k
%+\sum_{k=1}^{m_u} W_{n,k}^{(\beta-\alpha,\sigma)} U_k+\widehat{B}_n^\sigma U_0\right]\\
%&+Ah^{\beta}\left[\sum_{k=0}^n\omega_{n-k}^{(\beta)}U_k
%+\sum_{k=1}^{m_u} W_{n,k}^{(\beta,\sigma)} U_k+B_n^\sigma U_0\right]\\
%&+h^{\beta}\left[\sum_{k=0}^{n-1}\omega_{n-k}^{(\beta)}F_k
%+\sum_{k=1}^{m_f} W_{n,k}^{(\beta,\delta)} F_k+B_n^\delta F_0\right]\\
%&+h^{\beta}\omega_0^{(\beta)}\bigg[2F_{n-1}-F_{n-2}
%+\sum_{k=1}^{\widetilde{m}_f} \widehat{W}_{n,k}^{(f)} (F_k-F_0)\bigg],
%\end{aligned}\end{equation}
%where
%$
%\widehat{B}_n^\sigma=\frac{n^{\beta-\alpha}}{\Gamma(1+\beta-\alpha)}- \sum_{k=0}^n\omega_{n-k}^{(\beta-\alpha)} -\sum_{k=1}^{m_u}{W}_{n,k}^{(\beta-\alpha,\sigma)},
%$
%and the other terms in \eqref{scheme-ts-e-mul} are defined similarly as those in the IMEX-E scheme \eqref{scheme-ts-e}.
%Since the IMEX-T and IMEX-Trap schemes can be also derived readily, we do not present them  here.

{When using uniform stepsize, the  TS-I and TS-III schemes in \cite{CaoZhangK15} can be readily rewritten as the IMEX forms. For comparison,  we present  IMEX forms of the TS-I and TS-III schemes for \eqref{eq:multifde}.}

\textbf{IMEX form of TS-I:}
\begin{eqnarray*}
U_n&=&U_0+w^{(\beta-\alpha)}_{n,n}n^{\beta-\alpha}U_0-\sum_{j=1}^{n-1}w_{n,j}^{(\beta-\alpha)}
U_{j-\frac{1}{2}}-w^{(\beta-\alpha)}_{n,n}U_n\\
&&+A\left(\sum_{j=1}^{n-1}w_{n,j}^{(\beta)}U_{j-\frac{1}{2}}+w^{(\beta)}_{n,n}U_n\right)
+\sum_{j=1}^{n-1}w_{n,j}^{(\beta)}f
\left(t_{j-\frac{1}{2}},U_{j-\frac{1}{2}}\right)+w^{(\beta)}_{n,n}f(t_{n-1},U_{n-1}),
\end{eqnarray*}
where $U_{j-\frac{1}{2}}=\frac{U_{j-1}+U_{j}}{2}$ and
$w^{(\nu)}_{n,j}=\frac{h^{\nu}}{\Gamma(1+\nu)}[(n-j+1)^{\nu}-(n-j)^{\nu}]$, $\nu=\beta,\beta-\alpha$.

\textbf{IMEX form of TS-III:}
\begin{eqnarray*}
&&\sum_{j=1}^{n}b^{(\alpha)}_{n,j}
\left(U_{j}-U_{j-1}\right)+
\sum_{j=1}^{n}b^{(\beta)}_{n,j}\left(U_{j}-U_{j-1}\right)\\
&=&AU_n+f(t_{n},U_{n-1})+diag[f_u(t_n,U_{n-1})](U_n-U_{n-1}),
\end{eqnarray*}
where $b^{(\nu)}_{n,j}=\frac{h^{-\nu}}{\Gamma(2-\nu)}[(n-j+1)^{1-\nu}-(n-j)^{1-\nu}]$,
$\nu=\alpha,\beta$, and $diag[\cdot]$ denotes the diagonal matrix where the $k$-th diagonal element is the $k$-th element of a vector, $k=1,2,\cdots,m$.

We will test the schemes IMEX-E, IMEX-T and IMEX-E-Trap  in Example \ref{exm:3} and compare the numerical results with the TS-I, TS-III schemes and the predictor-corrector scheme \cite{DieFF02}. Convergence analysis will not be presented in this paper. However, it will be shown numerically that the IMEX-E, IMEX-T and IMEX-E-Trap schemes are of uniformly second-order accuracy for solving \eqref{eq:fdemultiterm}   with non-smooth solutions or  smooth solutions.

%%%%%%%%%%%%%%%%%%%%%%%%%%%%%%%%%%%%%%%%%%%%%%%%%%%%%%%%%%%

\section{Numerical examples}\label{sec:num}
We denote by  $U_n$ a numerical solution of the numerical methods in the present work with a time step size $h$ at $t_n=nh$, and we measure the errors in the following sense:
\[E^r_\infty(h)=\frac{\max_{0\leq n\leq  N}|U^{\rm ref}_n-U_n|}
{\max_{0\leq n\leq N}|U^{\rm ref}_n|},{\quad}
E^r_N(h)=\frac{|U^{\rm ref}_N-U_N|}{|U^{\rm ref}_N|}.\]
If the exact solution $u(t)$ is available, then we take $U^{\rm ref}_n=u(t_n)$; otherwise we obtain the reference solution  $U^{\rm ref}_n$ with the step size $h=2^{-15}$.

We will test accuracy and convergence order of the  IMEX-E and IMEX-T schemes for a stiff system and a nonlinear equation with non-smooth solutions in Examples \ref{exm:linearsys1} and  \ref{exm:nonlinear1}, respectively. In Example \ref{exm:3}, we will show the behavior of the IMEX-E, IMEX-T and IMEX-Trap schemes when solving a multi-term nonlinear stiff fractional differential system.  We will also compare our methods with the TS-I, TS-III schemes proposed in \cite{CaoZhangK15} (see their IMEX forms in Section \ref{sec4}), and the predictor-corrector (PC) scheme  developed in \cite{DieFF02} for $(\,^CD_{0}^{\beta}u)(t)=f(t,u),\; t\in(0,T],\; u(0)=u_0$:
\begin{eqnarray}
u^p_{n}&=&u_0+\sum_{j=0}^{n-1} w_{n-1,j}^{(\beta)} f(t_j,u_j),\label{eq:pc1}\\
u_{n}&=&u_0+h^{\beta}b_{n,n}^{(\beta)}f(t_{n},u_{n}^p)
+h^{\beta}\sum_{j=0}^{n-1}b_{n,j}^{(\beta)}f(t_j,u_j),\label{eq:pc2}
\end{eqnarray}
where  $w_{n,j}^{(\beta)}=\frac{h^{\beta}}{\Gamma(1+\beta)}[(n-j+1)^\beta-(n-j)^\beta],$ and
$b_{n,j}^{(\beta)}$ is defined in \eqref{eq:trapapp}.
%$$
%a_{n+1,j}=\frac{h^\beta}{\Gamma(\beta+2)}\left\{\begin{array}{ll}
%n^{\beta+1}-(n-\beta)(n+1)^\beta, & j=0; \\
%(n-j+2)^{\beta+1}+(n-j)^{\beta+1}-2(n-j+1)^{\beta+1}, & 1\leq j\leq n.
%\end{array}\right.
%$$

%, and consider the substepping based on IMEX-E in Example \ref{exm:substep1} to show flexibility of the time-splitting scheme. Then we compare the performance of the IMEX-E and IMEX-T scheme with the TS-I, TS-III scheme proposed in \cite{CaoZhangK15} for solving the TFDEs with smooth solutions in Example \ref{exm:3}.

\vspace{0.3em}
\begin{exm}[Stiff fractional ordinary differential system]\label{exm:linearsys1}
%%%%%%%%%%%%%%%%%%%%%%%%%%%%%%%%%%
\begin{eqnarray}\label{eq:linearsys1}
(\,^CD_0^{\beta}u)(t)&=&{Au(t)+Bu(t)}+g(t),\; t\in(0,T],\\
u(0)&=&u_0.\nonumber
\end{eqnarray}
\end{exm}
In this example, we take $f(t,u)=Bu(t)+g(t)$ and use the IMEX-E scheme \eqref{scheme-ts-e} to solve this system.
Take $u_0=(u_{01},u_{02},u_{03})^\top= (1,1,1)^{\top} $,
\begin{equation}\label{eq:stiff}
A=\left(
\begin{array}{ccc}
-10000 & 0 & 1\\
-0.05 & -0.08&-0.2 \\
1  &  0  &-1\\
\end{array}
\right),\;\;B=\left(
\begin{array}{ccc}
-0.6 & 0 & 0.2 \\
-0.1 & -0.2 & 0\\
0  & -0.5 & -0.8\\
\end{array}
\right),
\end{equation}
and
\begin{equation}\label{eq:stiff_f}
g(t)=\left(
\begin{array}{c}
a_1\Gamma_1t^{{\sigma}_1-\beta}+a_2\Gamma_2t^{{\sigma}_2-\beta}\\
a_3\Gamma_3t^{{\sigma}_3-\beta}+ a_4\Gamma_4t^{{\sigma}_4-\beta}\\
a_5\Gamma_5t^{{\sigma}_5-\beta}+a_6\Gamma_6t^{{\sigma}_6-\beta}\\
\end{array}\right)-(A+B)\left(
\begin{array}{c}
a_1t^{{\sigma}_1}+a_2t^{{\sigma}_2}+u_{01}\\
a_3t^{{\sigma}_3}+a_4t^{{\sigma}_4}+u_{02}\\
a_5t^{{\sigma}_5}+a_6t^{{\sigma}_6}+u_{03}\\
\end{array}
\right),
\end{equation}
where $\Gamma_k=\frac{\Gamma({\sigma}_k+1)}{\Gamma({\sigma}_k+1-\beta)}\,(1\leq k \leq 6)$.
Then  the exact solution of \eqref{eq:linearsys1} is
\begin{equation}\label{eq:stiff_ext}
u(t)= ( a_1t^{{\sigma}_1}+a_2t^{{\sigma}_2}+u_{01},   a_3t^{{\sigma}_3}+a_4t^{{\sigma}_4}+u_{02},    a_5t^{{\sigma}_5}+a_6t^{{\sigma}_6}+u_{03})^\top.
\end{equation}

We take ${\sigma}_1=\beta,\;{\sigma}_2=2\beta,\;{\sigma}_3=1+\beta,\; {\sigma}_4=5\beta,\; {\sigma}_5=2,\; {\sigma}_6=2+\beta,$ and
$a_1=0.5, \;a_2=0.8,\;a_3=1,\; a_4=1,\;a_5=1,\;a_6=1$ in the numerical computation.

\begin{table}[!htb]
\centering
\caption{
Relative error and convergence rate of the IMEX-E scheme \eqref{scheme-ts-e} with different correction terms for the stiff system \eqref{eq:linearsys1} (Example \ref{exm:linearsys1}), $T=1$, $\beta=0.1$.}\label{tbl:linsys1-1}
\scalebox{0.8}{
\begin{tabular}{c|cc|cc|cc|cc|cc}
\hline
\hline
\multirow{2}{*}{$h$}& \multicolumn{2}{|c|}{$m=0$} & \multicolumn{2}{|c}{$m=1$}& \multicolumn{2}{|c|}{$m=2$} & \multicolumn{2}{|c}{$m=3$} & \multicolumn{2}{|c}{$m=4$} \\

\cline{2-3}
\cline{4-5}
\cline{6-7}
\cline{8-9}
\cline{10-11}

~& $ E^r_\infty(h)$ &order  &  $ E^r_\infty(h)$ &order   & $ E^r_\infty(h)$ &order  & $ E^r_\infty(h)$ &order& $ E^r_\infty(h)$ &order   \\
\hline

 $2^{-10}$   & 9.45e-3	 &  0.24 & 1.13e-3 &  0.48 & 1.33e-4 & 0.61 & 2.43e-6 & 1.19  & 2.27e-7 & 2.05  \\
 $2^{-11}$   & 8.02e-3	 &  0.23 & 8.10e-4 &  0.46 & 8.68e-5 & 0.60 & 1.07e-6 & 1.17  & 5.46e-8 & 2.05    \\
 $2^{-12}$   & 6.83e-3	 &  0.23 & 5.89e-4 &  0.43 & 5.74e-5 & 0.58 & 4.73e-7 & 1.16  & 1.32e-8 & 2.05   \\
 $2^{-13}$   & 5.83e-3	 &  *    & 4.36e-4 &  *    & 3.83e-5 &   *  & 2.11e-7 & *     & 3.17e-9 & *       \\
\hline
\end{tabular}}
\end{table}

\begin{table}[!htb]
\centering
\caption{
Relative error and convergence rate of the IMEX-E scheme \eqref{scheme-ts-e} with different correction terms for the stiff system \eqref{eq:linearsys1} (Example \ref{exm:linearsys1}), $T=1$, $\beta=0.5$.}\label{tbl:linsys1-2}
\scalebox{0.85}{
\begin{tabular}{c|cc|cc|cc|cc}
\hline
	\hline
	\multirow{2}{*}{$h$}& \multicolumn{2}{|c|}{$m=0$} & \multicolumn{2}{|c}{$m=1$}& \multicolumn{2}{|c|}{$m=2$} & \multicolumn{2}{|c}{$m=3$}\\
	\cline{2-3}
	\cline{4-5}
	\cline{6-7}
	\cline{8-9}		
~& $ E^r_\infty(h)$ &order  &  $ E^r_\infty(h)$ &order   & $ E^r_\infty(h)$ &order  & $ E^r_\infty(h)$ &order \\
\hline		
$2^{-10}$   & 8.27e-4	 &  0.50 & 3.71e-5 &  1.00 & 1.06e-7 & 2.07 & 2.95e-8	 & 1.80   \\
$2^{-11}$   & 5.84e-4	 &  0.50 & 1.86e-5 &  1.00 & 2.52e-8 & 2.05 & 8.46e-9	 & 1.87     \\		$2^{-12}$   & 4.12e-4	 &  0.50 & 9.28e-6 &  1.00 & 6.11e-9 & 2.03 & 2.32e-9	 & 1.91    \\		$2^{-13}$   & 2.91e-4	 &  *    & 4.63e-6 &  *    & 1.49e-9 &   *  & 6.17e-10   & *        \\
\hline
\end{tabular}}
\end{table}

	\begin{table}[!htb]
	\centering
{	\caption{Condition numbers and residuals of the linear system \eqref{eq:van1} (Example \ref{exm:linearsys1}),  ${\sigma}_1=\beta,\;{\sigma}_2=2\beta,\;{\sigma}_3=1+\beta,\; {\sigma}_4=5\beta$. }\label{tbl:con1}
	\scalebox{0.85}{
		\begin{tabular}{|c|c|c|c|c|}
			\hline
			 	$\beta$ & - & $m=2$ & $m=3$ & $m=4$\\
		\hline 
	\multirow{2}{*}{0.1}& condition numbers& 6.20e+01& 1.70e+03& 2.85e+04\\
	&residuals & 6.94e-18 & 5.55e-17& 1.11e-15\\\hline 
	\multirow{2}{*}{0.5} & condition numbers& 1.74e+01 & 2.65e+02& 5.11e+03 \\
	& residuals& 5.54e-17 & 1.10e-16& 8.53e-14\\
			\hline
		\end{tabular}}}
	\end{table}

\begin{figure}[!h]
		\caption{Comparison of relative errors for the IMEX-E  scheme with different correction terms for the stiff    system \eqref{eq:linearsys1}  (Example \ref{exm:linearsys1}), $h=2^{-12},\;T=1$. All needed starting values are given in advance.}
		\label{fig:stiffsys1}
	\centering
	\epsfig{figure=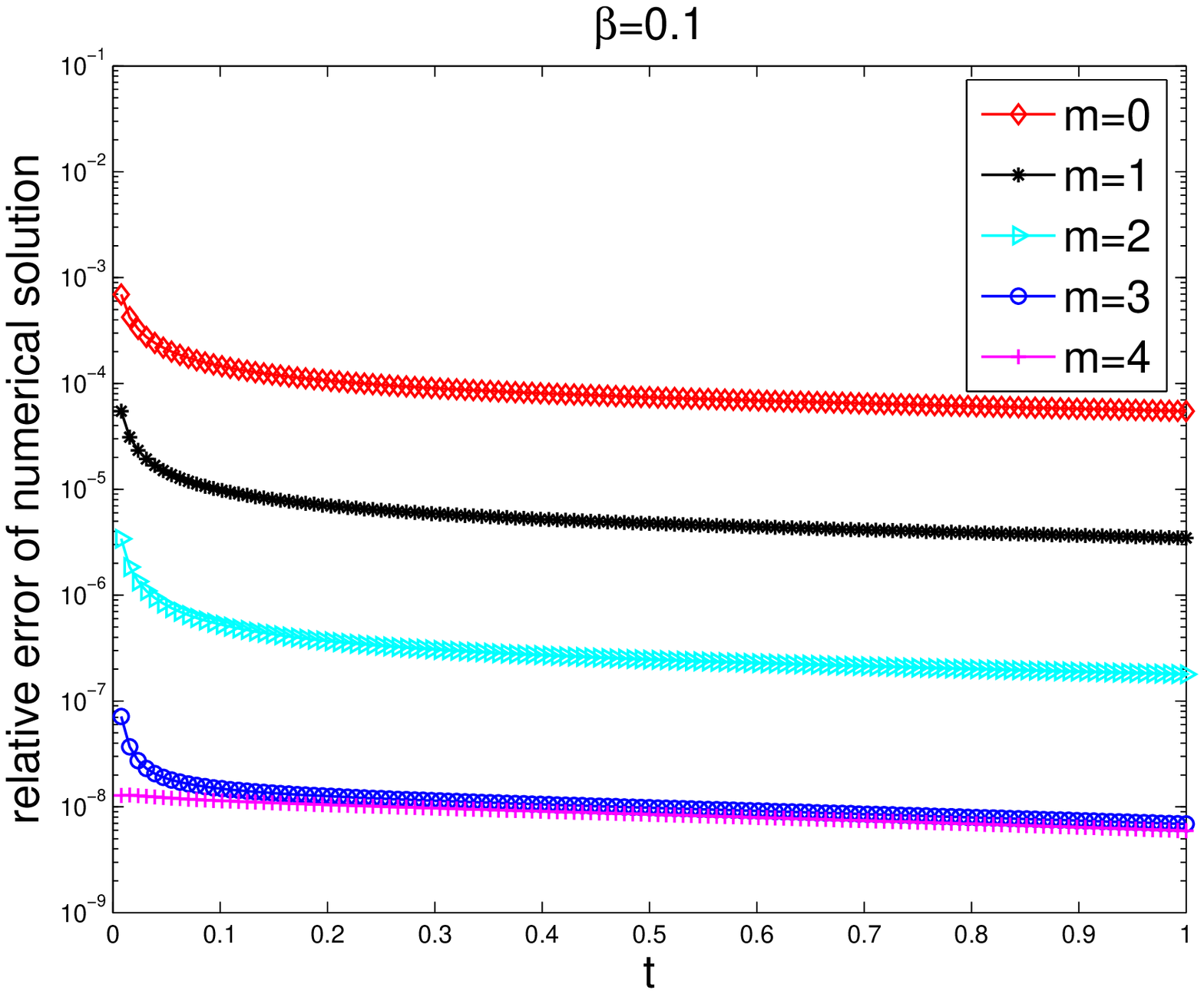,width=6cm}~\epsfig{figure=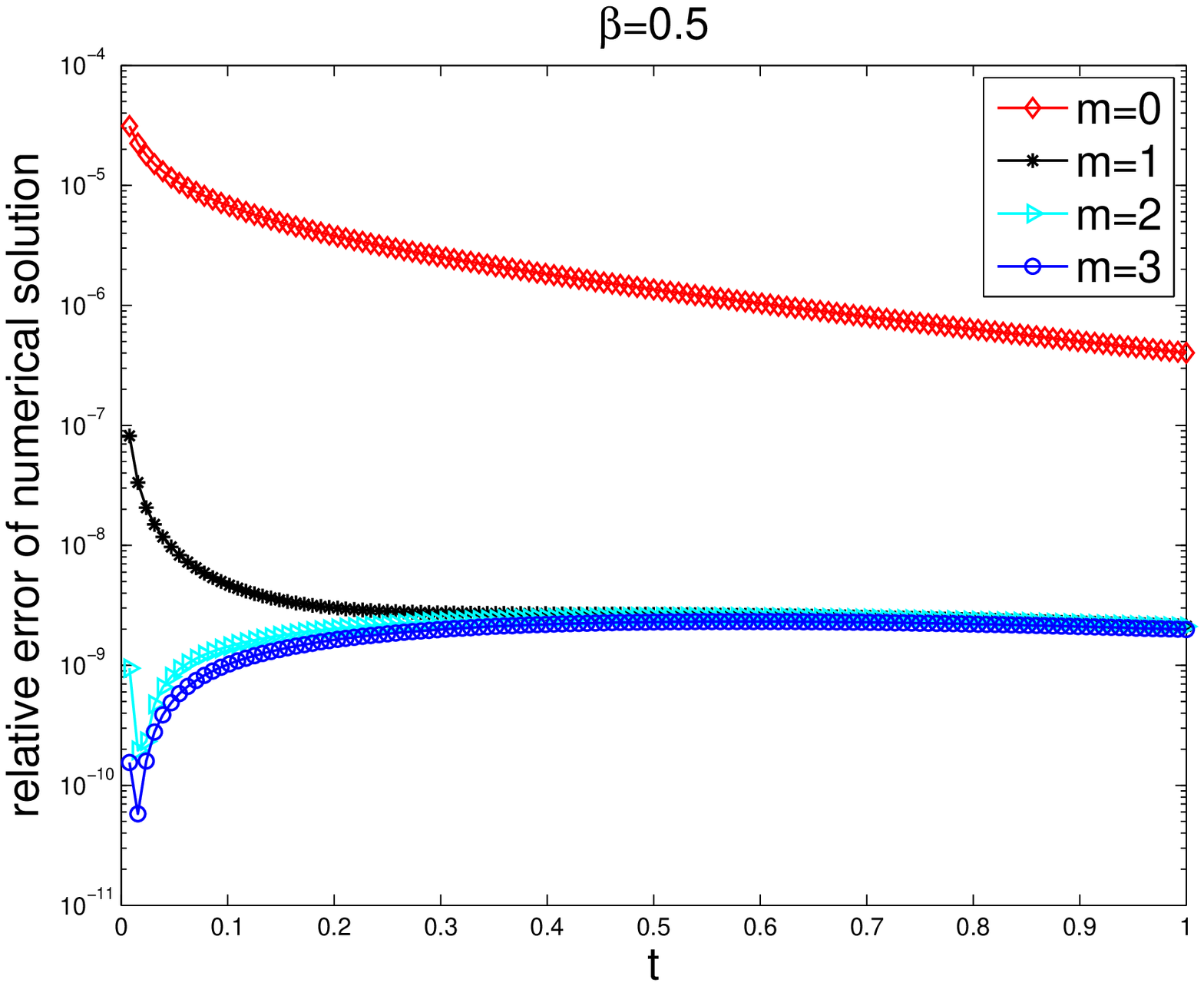,width=6cm}
\end{figure}

Tables \ref{tbl:linsys1-1} and \ref{tbl:linsys1-2} show that the use of starting values in the  IMEX-E scheme is crucial for both reducing relative errors and getting uniformly second-order accuracy.  We observe that the fractional system \eqref{eq:linearsys1} with smaller fractional order $\beta=0.1$ requires more correction terms than that with bigger fractional order $\beta=0.5$, which is consistent with Theorem \ref{thm:convergence}. Moreover, from the data in the last column of Table \ref{tbl:linsys1-2}, we find that no gain of convergence rate is observed when the number of correction terms exceeds some optimal number that can be derived from  Theorem \ref{thm:convergence}.
{ Specifically, the global second-order accuracy can be  obtained when $\min\{\sigma_{m+1}+\beta,\delta_{m+1}+\beta\}\geq 2$.
For $\beta=0.1$, the smallest $m$ to get the global second-order accuracy
is $4$, while for $\beta=0.5$, the optimal $m$ is $2$; numerical results in Tables
\ref{tbl:linsys1-1}--\ref{tbl:linsys1-2} verify the theoretical results in
Theorem \ref{thm:convergence}.}
%{In fact, by Theorem \ref{thm:convergence},  the IMEX-E scheme is of second-order convergence when $q=2$. To obtain the second-order accuracy, we accurately calculate the singular terms of solution, which cause $q<2$, by applying correction terms. It determines how many and which terms need to be corrected when solving the system \ref{eq:linearsys1} in Example \ref{exm:linearsys1}. For instance, when $\beta=0.5$, we do not need to apply correction terms to $\sigma_4=5\beta$, since $5\beta+\beta>2$, but when $\beta=0.1$, we do need to correct it.  }
Data in the last column of Table \ref{tbl:linsys1-2} also imply that too many correction terms are not always helpful, which can be also illustrated from Lemma \ref{lem:weight2}, where  the starting weights (see, e.g.  $W_{n,k}^{(\beta,\sigma)}$ in \eqref{eq:weight1}) will be large when $\sigma_r$ is sufficiently large.
{ However, from the condition numbers and residuals of \eqref{eq:van1} in Table \ref{tbl:con1},  we observe that relatively high accuracy of starting weights can be obtained even when the condition number is large if there are only a few correction terms.
Here and in Table \ref{tbl:con2}, the residual is computed by
$$\max_{1\leq r\leq m, 1\leq n\leq 100}\left|\sum_{k=1}^{m} W_{n,k}^{(\beta,u)}k^{\sigma_r}-\frac{\Gamma(\sigma_r+1)}{\Gamma(\sigma_r+1+\beta)}n^{\sigma_r+\beta}
 +\sum_{k=0}^n\omega^{(\beta)}_{n-k}k^{\sigma_r}\right|.$$}
%how to choose $t^{\sigma_r}$ which needs to be corrected, and too large $\sigma_r$ may lead to large weights of the correction terms.

%We also solve \eqref{eq:linearsys1}--\eqref{eq:stiff} by the predictor-corrector scheme proposed in \cite{DieFF04}£º
%\begin{eqnarray}
%u_{n+1}&=&u_0+\frac{h^\beta}{\Gamma(\beta+2)}\left(\sum_{j=0}^na_{j,n+1}u(t_j,u_j)+a_{n+1,n+1}u(t_{n+1},u^P_{n+1})\right),\label{eq:adams1}\\
%u_{n+1}^P&=&u_0+\frac{1}{\Gamma(\beta+1)}\sum_{j=0}^nb_{j,n+1}u(t_j,u_j),\label{eq:adams2}
%\end{eqnarray}
%where
%\begin{eqnarray*}
%a_{j,n+1}&=&
%\left\{\begin{array}{ll}
%	n^{\beta+1}-(n-\beta)(n+1)^\beta,& if\,j=0,\\
%(n-j+2)^{\beta+1}+(n-j)^{\beta+1}-2(n-j+1)^{\beta+1},& if\,1\leq j\leq n,\\
%1, & if\,j=n+1
%	\end{array}\right.\\
%b_{j,n+1}&=&(n+1-j)^\beta-(n-j)^\beta.
%\end{eqnarray*}

We further solve \eqref{eq:linearsys1} by the predictor-corrector scheme \eqref{eq:pc1}-\eqref{eq:pc2} for comparison.
In our tests, numerical solutions of the predictor-corrector scheme  blow up very quickly for $\beta=0.1$ and $\beta=0.5$, even if we take very small step size $h=2^{-16}$. For  bigger $\beta$ (i.e.,  $\beta=0.95$), the predictor-corrector scheme cannot work either, except that we take very small stepsize (i.e.,  $h=2^{-13}$). We do not present all these results here.

Figure \ref{fig:stiffsys1} shows the asymptotic behavior of relative errors for the IMEX-E scheme with different number of correction terms. It is shown that suitable correction terms can improve the accuracy greatly. Moreover, there exists an optimal number of correction terms for both $\beta=0.1$ and $\beta=0.5$. If one applies more correction terms than the optimal choice, the accuracy may not be further improved. In this example, all necessary starting values have been given in advance
(We have analytical solutions). However, in practice, we have to calculate starting values numerically using high-order methods or using small stepsizes, which will be shown in the following
example.

%will take different  computational time when using different numbers of correction terms. We will show it in Example \ref{exm:nonlinear1}.

%%%%%%%%%%%%%%%%%%%%%%%%%%%%%%%%%%
\vspace{0.3em}
\begin{exm}[Nonlinear fractional ordinary differential equation]\label{exm:nonlinear1}
	\begin{eqnarray}\label{eq:nonlinear1-1}
	(\,^CD_0^{\beta}u)(t)&=&{\lambda u(t)}+\rho u(1-u^2)+g(t), t\in(0,T], \quad 	u(0) = u_0.
	\end{eqnarray}
\end{exm}
\begin{itemize}
  \item Case I:  Take $\lambda=-3$, $\rho=0.8$. Choose suitable
  $g(t)$ such that the  solution to \eqref{eq:nonlinear1-1} is $\displaystyle u(t) = u_0+\sum_{k=1}^6t^{\sigma_k},$
  where $u_0=2$, $\sigma_k=k\beta\,(1\leq k \leq 5)$ and $\sigma_6=2+\beta$.
  \item Case II:  Take $g(t)=0$, $\lambda=-3$, $\rho=0.8$, and $u_0=2$.
  From Lemma \ref{lem:diethm1}, we know that the analytical solution $u(t)$ satisfies \eqref{eq:structure}. So $\sigma_k$ in \eqref{eq:analsol} satisfies
  $\sigma_k\in\{i+j\beta,i=0,1,...,j=1,2,...\}$.
  \end{itemize}
  As we do not have an exact solution for Case II,  we calculate a reference solution by the considered scheme with very small step size $h=2^{-15}$.

Table \ref{tbl:fzeng1} shows the maximum relative errors of the IMEX-T scheme \eqref{scheme-ts-t} with different time stepsizes for solving the nonlinear equation \eqref{eq:nonlinear1-1} (Case I) when $\beta=0.15$.  { We also apply the IMEX-E scheme \eqref{scheme-ts-e} to  this problem, and observe that numerical solutions blow up even for the small step size   $h=2^{-14}$ (results not presented here).}

\begin{figure}[!t]
		\caption{Asymptotic relative error for the IMEX-E scheme \eqref{scheme-ts-e} and IMEX-T scheme \eqref{scheme-ts-t} with different correction terms for the nonlinear  equation \eqref{eq:nonlinear1-1} (Example \ref{exm:nonlinear1}, Case II), $T=8,\;\beta=0.15$.}
		\label{fig:nl-2}
	\centering
	\epsfig{figure=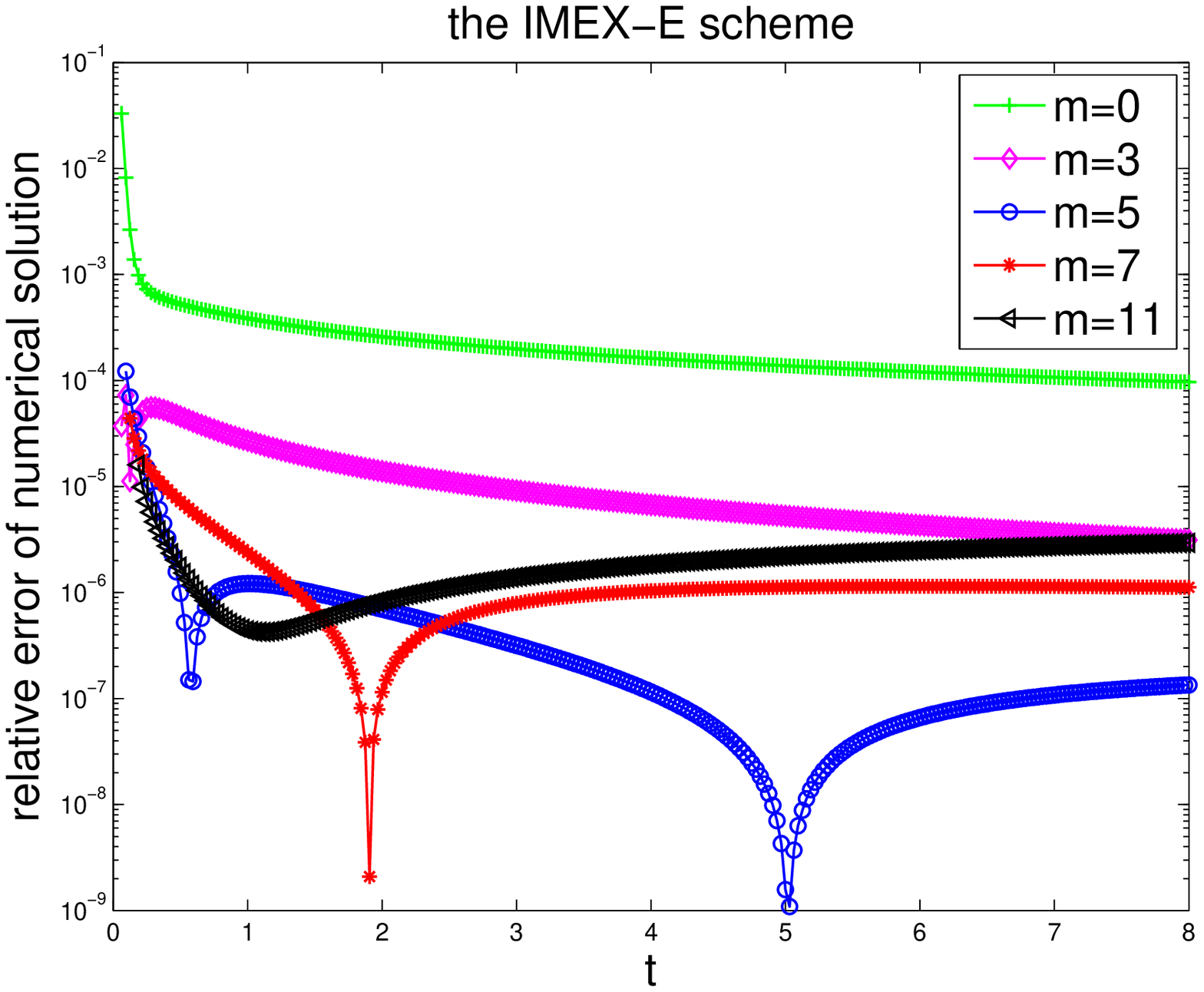,width=6cm}~\epsfig{figure=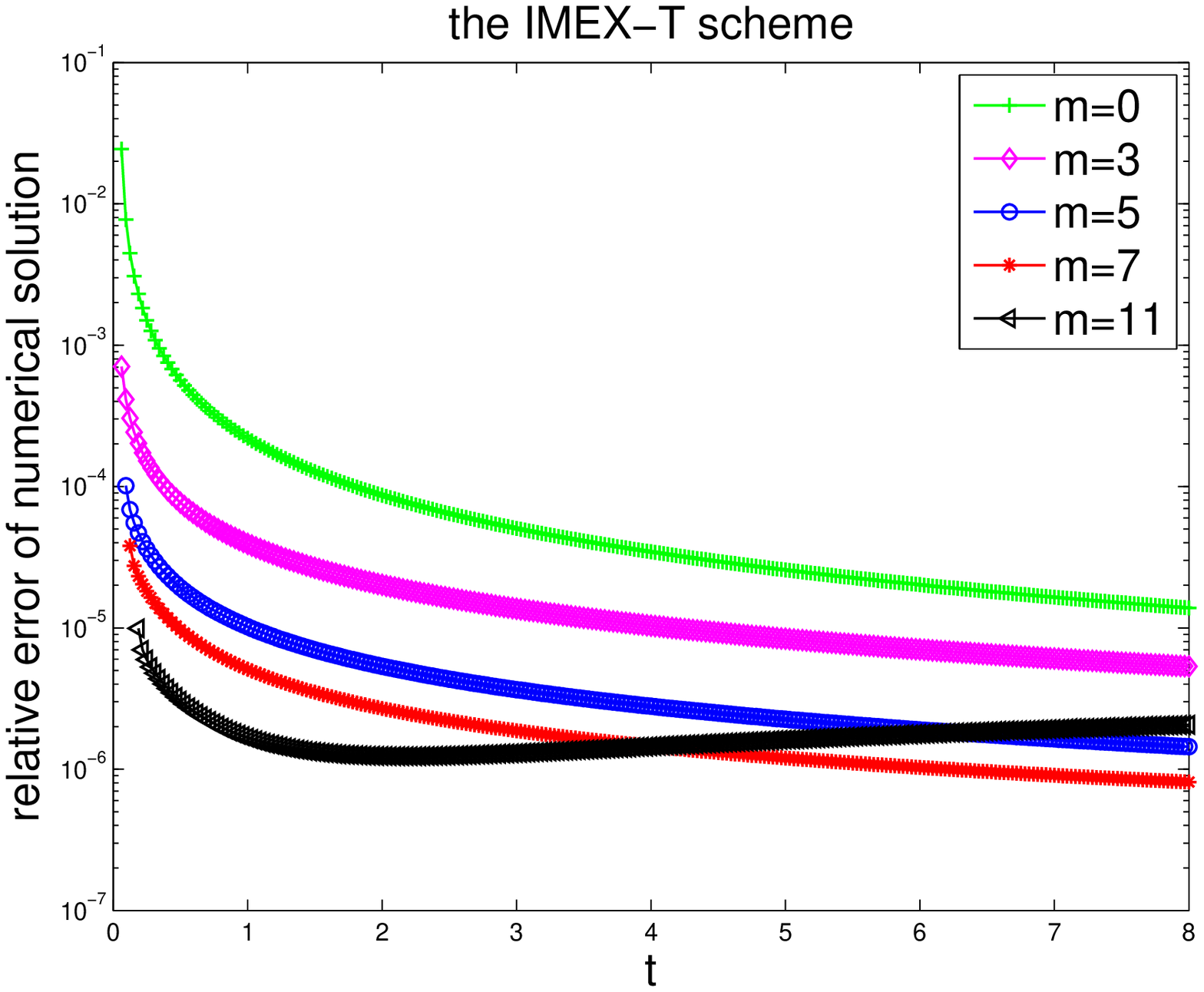,width=6cm}
\end{figure}

Table \ref{tbl:nlinear-2} and Figure \ref{fig:nl-2}  show the relative errors and convergence rates of the IMEX-E and IMEX-T scheme for solving the nonlinear equation \eqref{eq:nonlinear1-1} (Case II). It is shown that we can apply proper correction terms for both schemes to improve their convergence order up to second-order though we do not know what the exact solution is. {Here we use the adaptive step size $h^2$ to calculate starting values of correction terms in our schemes for different stepsize $h$.} Numerical tests show that it works well. {In Table \ref{tbl:con2}, we also present the condition numbers and residuals of system \eqref{eq:van1} for  different number of correction terms. }% in Table \ref{tbl:nlinear-2}. }

When using the predictor-corrector scheme  \eqref{eq:pc1}-\eqref{eq:pc2} to solve \eqref{eq:nonlinear1-1} (Case II) in this example, the numerical solution blows up for $\beta=0.15$ and is of convergence order $\beta$ for $\beta=0.95$; see Table \ref{tbl:adams}.

\begin{table}[!tb]
\centering
\caption{
Relative error and convergence rate of the IMEX-T scheme \eqref{scheme-ts-t} for  the nonlinear equation \eqref{eq:nonlinear1-1} (Example \ref{exm:nonlinear1}, Case I), $T=8$, $\beta=0.15$.}\label{tbl:fzeng1}
\scalebox{0.85}{
\begin{tabular}{c|cc|cc|cc|cc|cc}
\hline
	\hline
	\multirow{2}{*}{$h$}& \multicolumn{2}{|c|}{$m=0$} & \multicolumn{2}{|c}{$m=1$}& \multicolumn{2}{|c|}{$m=2$} & \multicolumn{2}{|c}{$m=3$}
& \multicolumn{2}{|c}{$m=4$}\\
	\cline{2-3}	\cline{4-5}	\cline{6-7}	\cline{8-9}	\cline{10-11}		
~& $ E^r_\infty(h)$ &order  &  $ E^r_\infty(h)$ &order   & $ E^r_\infty(h)$ &order
& $ E^r_\infty(h)$ &order& $ E^r_\infty(h)$ &order \\
\hline		
$2^{-5}$&8.99e-4&0.11   &6.05e-4&0.94&1.07e-3&1.71&6.20e-4&2.07&6.04e-4&2.04   \\
$2^{-6}$&9.70e-4&0.30&3.14e-4&0.08&3.27e-4&1.36&1.46e-4&2.01&1.46e-4&2.10     \\		
$2^{-7}$&7.87e-4&0.38&2.95e-4&0.37&1.27e-4&1.02&3.64e-5&1.77&3.40e-5&2.13    \\	
$2^{-8}$&6.03e-4&0.39&2.28e-4&0.43&6.26e-5&0.81&1.06e-5&1.45&7.76e-6&2.14     \\
$2^{-9}$&4.59e-4& *  &1.68e-4& *  &3.56e-5&*   &3.88e-6&*   &1.76e-6&*       \\
\hline
\end{tabular}}
\end{table}

\begin{table}[!h]
	\centering
	\caption{
		Relative errors and convergence rate of the IMEX-E scheme \eqref{scheme-ts-e} and the IMEX-T scheme \eqref{scheme-ts-t}. Comparison of the schemes  with different correction terms for the nonlinear equation \eqref{eq:nonlinear1-1} (Example \ref{exm:nonlinear1}, Case II), $T=8$, $\beta=0.15$.}\label{tbl:nlinear-2}
	\scalebox{0.85}{
		\begin{tabular}{c|c|ccc|ccc}
			\hline
			\hline
			\multirow{2}{*}{$m$}& \multirow{2}{*}{$h$}& \multicolumn{3}{|c|}{IMEX-E} & \multicolumn{3}{|c}{IMEX-T}\\
			\cline{3-5}
			\cline{6-8}
			~&~& $ E^r_\infty(h)$ & order & cputime &  $ E^r_\infty(h)$ & order  & cpu time  \\\hline	
			$0$ &	$2^{-4}$   &2.82e-01&	0.10	&	0.02 & 1.51e-01 &	0.21	&	0.03	\\
			&	$2^{-5}$       &2.63e-01&	0.10	&	0.08 & 1.31e-01 &	0.21	&	0.05	 \\	
			&	$2^{-6}$       &2.46e-01&	0.09	&	0.13 & 1.13e-01 &	0.21	&	0.13	\\
			&	$2^{-7}$       &2.31e-01&	*   	&	0.27 & 9.82e-02 &	*   	&	0.28     \\\hline
			$3$ &	$2^{-4}$  & 1.28e-03&	1.80	&	0.03 & 2.02e-03	&	1.01	&	0.03	\\		
			&	$2^{-5}$      & 3.68e-04&	0.94	&	0.05 & 1.00e-03	&	0.51	&	0.08	 \\	
			&	$2^{-6}$      & 1.92e-04&   0.16	&	0.12 & 7.06e-04	&	0.37	&	0.14	\\	
			&	$2^{-7}$      & 1.72e-04&	*   	&	0.29 & 5.48e-04	&	*   	&	0.29     \\\hline
			$5$&$2^{-4}$      &8.23e-04&	1.26	&	0.03 & 8.52e-04	&	1.68	&	0.03	\\		
			&	$2^{-5}$      &3.45e-04	&	0.93	&	0.08 & 2.67e-04	&	1.39	&	0.07	 \\		
			&	$2^{-6}$      &1.81e-04	&	0.60	&	0.14 & 1.01e-04	&	1.05	&	0.16	\\		
			&	$2^{-7}$      &1.20e-04	&	*   	&	0.35 & 4.91e-05	&  *    	&	0.36     \\		\hline
			$7$&$2^{-4}$      &6.07e-04&	1.58	&	0.06 & 5.83e-04&	1.79&	0.04	\\		
			&	$2^{-5}$      &2.03e-04&    1.60	&	0.07 &1.68e-04	&2.14	&	0.08	 \\		
			&	$2^{-6}$      &6.71e-05&    1.40	&	0.16 &3.81e-05	&2.74	&	0.16	\\		
			&	$2^{-7}$      &2.55e-05&    *   	&	0.44 &5.70e-06	&  *    	&	0.38     \\		\hline
			$11$&$2^{-4}$     &3.57e-04&	1.96	&	0.05 & 3.55e-04	&	2.05	&	0.06	\\		
			&	$2^{-5}$      &9.20e-05	&	1.86	&	0.09 & 8.56e-05	&	2.45	&	0.09	 \\		
			&	$2^{-6}$      &2.53e-05	&	1.99	&	0.24 & 1.57e-05	&	2.10	&	0.22	\\		
			&	$2^{-7}$      &6.36e-06	&	*   	&	0.55 & 3.66e-06	&	*   	&	0.64     \\		
			\hline
		\end{tabular}}
	\end{table}
\begin{table}[!htb]
	\centering
	{	\caption{Condition numbers and residuals of the linear system \eqref{eq:van1} (Example \ref{exm:nonlinear1}, Case II), $\beta=0.15$,  ${\sigma}_k=k\beta,\; k\geq 1$.}\label{tbl:con2}
		\scalebox{0.85}{
			\begin{tabular}{|c|c|c|c|c|}
				\hline	
		$m$ & $m=3$ & $m=5$ & $m=7$ & $m=11$\\
				\hline
				condition numbers& 2.06e+03& 3.32e+06& 6.43e+09& 2.54e+16\\\hline
				residuals & 5.55e-17 & 4.88e-15& 6.81e-13& 6.74e-08\\\hline
			\end{tabular}}}
		\end{table}
	
\begin{table}[!hbt]
\centering
\caption{Error  and convergence order of the predictor-corrector scheme \eqref{eq:pc1}--\eqref{eq:pc2} at $t=8$ for
the nonlinear   equation \eqref{eq:nonlinear1-1} (Example \ref{exm:nonlinear1}, Case II).}\label{tbl:adams}
\scalebox{0.85}{
\begin{tabular}{c|c|c|ccc|ccc}
\hline\hline
\multirow{2}{*}{$h$}& $\beta=0.15$  & \multirow{2}{*}{$h$}& \multicolumn{3}{|c|}{$\beta=0.55$}
& \multicolumn{3}{|c}{$\beta=0.95$}\\
\cline{2-2}\cline{4-6}\cline{7-9}
~& $ E^r_N(h)$& ~ &  $ E^r_N(h)$ & order & cpu time &  $ E^r_N(h)$ & order & cpu time  \\
\hline
$2^{-12}$ & Calculation Failed&	$2^{-5}$ &	1.62e-01&0.87   &0.05 & 2.50e-03 & 1.18& 0.05\\
$2^{-13}$ & Calculation Failed&	$2^{-6}$ &	8.86e-02&0.86   &0.13 & 1.10e-03 & 1.04& 0.14\\
$2^{-14}$ & Calculation Failed&	$2^{-7}$ &	4.87e-02&0.81   &0.38 & 5.38e-04 & 0.98& 0.39\\
$2^{-15}$ & Calculation Failed&	$2^{-8}$ &	2.78e-02&*      &2.04 & 2.72e-04 & *   & 2.09\\
\hline
\end{tabular}}

\end{table}

\vspace{0.3em}
\begin{exm}[A comparison of IMEX-T, IMEX-E, IMEX-E-Trap and TS-I, TS-III in \cite{CaoZhangK15}, PC scheme in \cite{DieFF02} for solving a multi-term nonlinear stiff fractional differential system]\label{exm:3}
\end{exm}
In this example, we solve the following multi-term nonlinear stiff fractional differential system:
\begin{eqnarray}\label{eq:fdemultiterm}
(\,^CD_{0}^{\alpha}u)(t)+(\,^CD_{0}^{\beta}u)(t)&=&A u(t)+B\sin u +g(t),\; t\in(0,T],\; u(0)=u_0,
\end{eqnarray}
where $\beta=0.55$, $\alpha=0.4$,
\begin{equation*}
A=\left(
\begin{array}{ccc}
-1000 & 100 \\
0 & -0.1
\end{array}
\right),\;\;B=\left(
\begin{array}{ccc}
1 & 0  \\
0 & 3 \\
\end{array}
\right).
\end{equation*}

\begin{figure}[!t]
			\caption{A comparison of asymptotic behavior of  relative errors of the IMEX-T, IMEX-E, IMEX-E-Trap, TS-I, TS-III, and PC schemes  for the multi-term nonlinear stiff system \eqref{eq:fdemultiterm} (Example \ref{exm:3}), $T=1,\,h=2^{-8}$.}
			\label{fig:comp4_1}
		\centering
		\includegraphics[scale=.32]{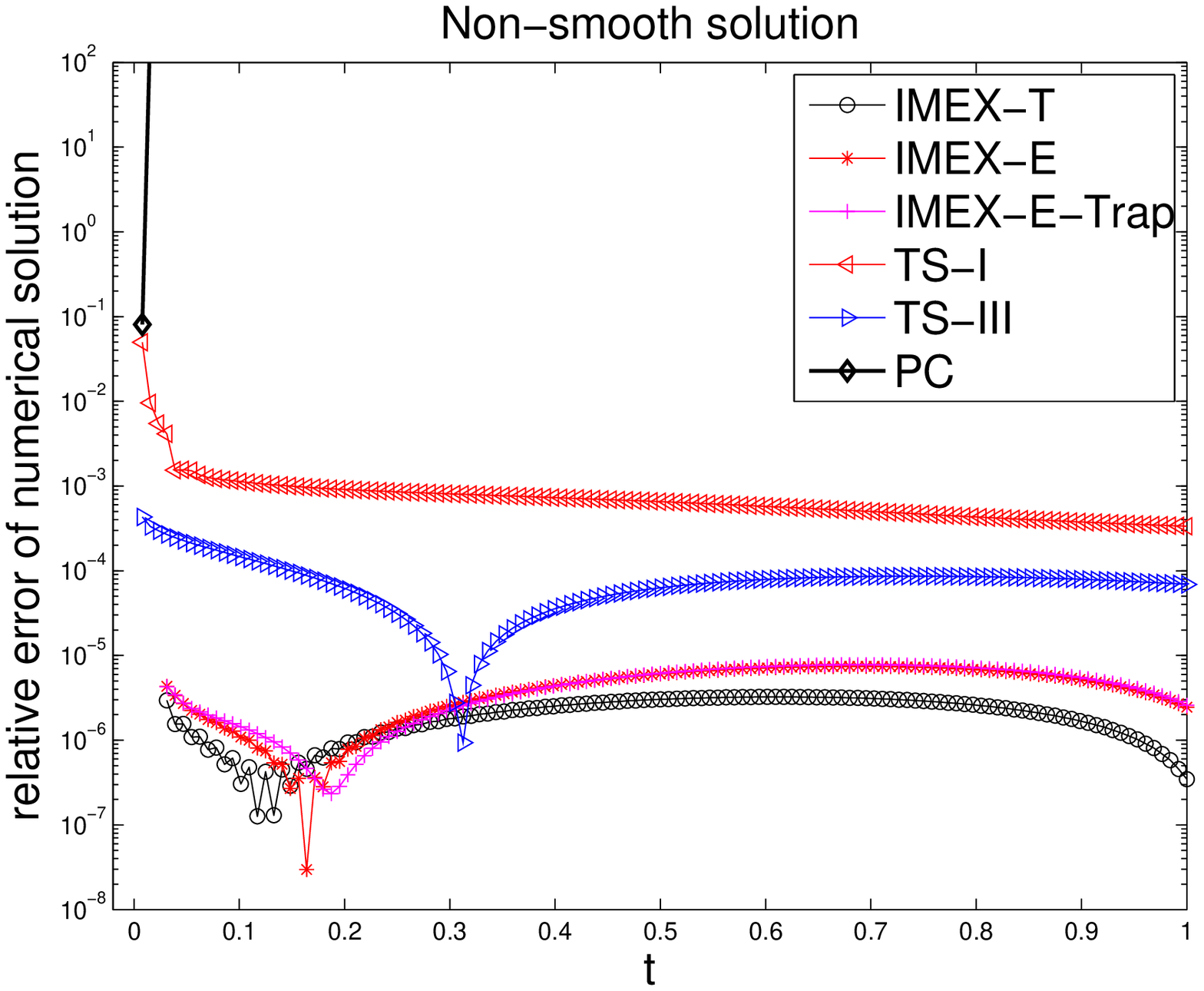}~\includegraphics[scale=.32]{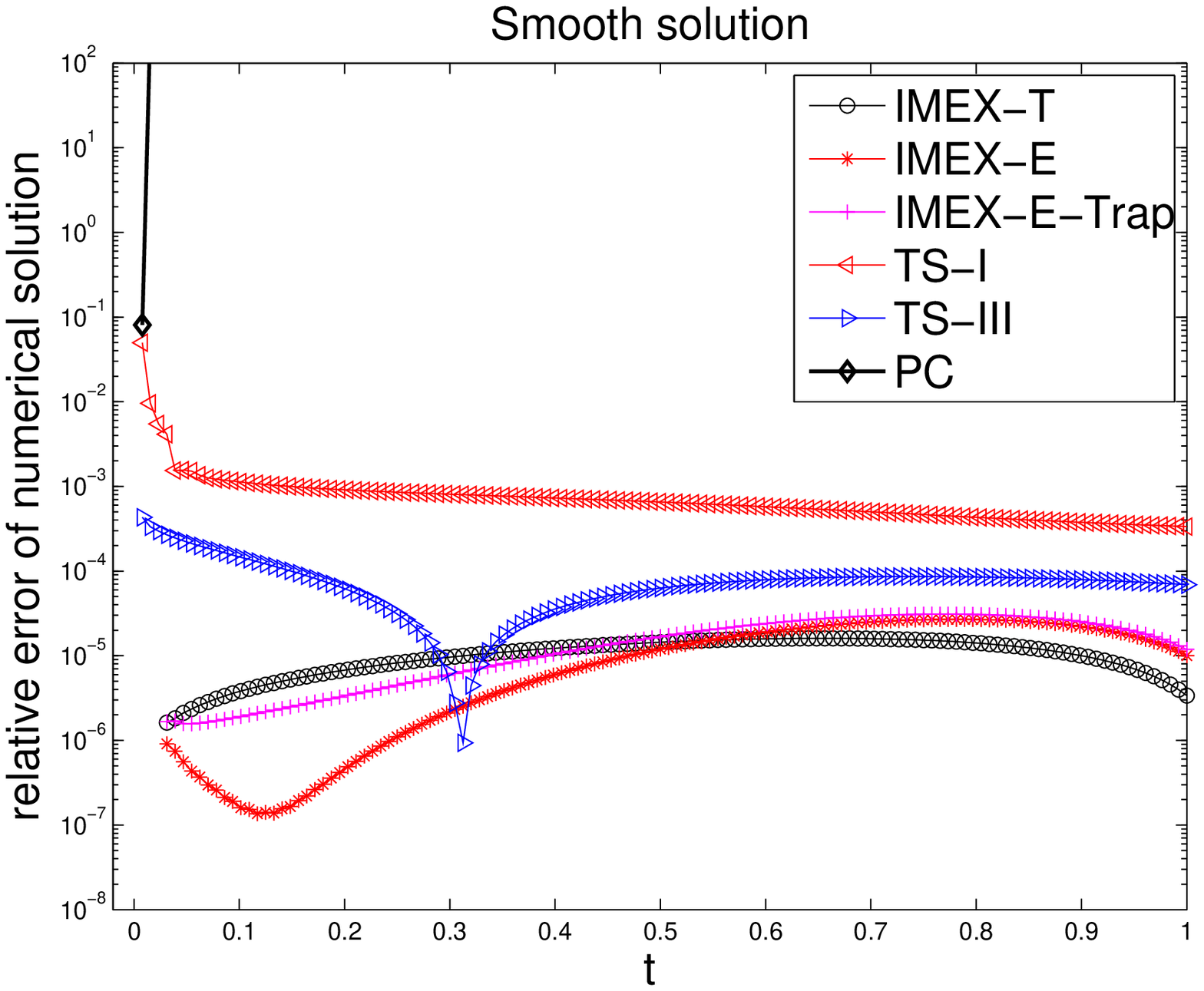}
	
	\end{figure}

We choose suitable $g(t)$ such that the above system has
\[ \text{
\textbf{Case I:}  }
 u(t)=\left[\begin{array}{c}
t^{0.55}+t^{1.15}+1\\
t^{0.9}+t^{2.55}+1
\end{array}
\right]; \quad
\text{
\textbf{Case II:}   }
 u(t)=\left[\begin{array}{c}
t^{2.55}+t^{2.15}+1\\
t^{3}+t^{2.55}+1
\end{array}
\right]. \]
Data in Table \ref{tbl:nlinear-comp4} show that the IMEX-T, IMEX-E and IMEX-E-Trap schemes work well for solving the problem \eqref{eq:fdemultiterm}, with both non-smooth solution and smooth solution, obtaining second-order accuracy numerical solutions. It is also shown that the TS-I and TS-III schemes solve the problem \eqref{eq:fdemultiterm} with smooth solution and get numerical solutions being of $1+\beta$ and $2-\beta$ order accuracy respectively, while for the non-smooth case, the accuracy is low and the convergence order of the TS-I and TS-III schemes is just $\beta$. It is worth to mention that in order to keep the uniformly second-order accuracy, we have applied suitable correction terms ($m_u=\widetilde{m}_u=3,\,m_f=\widetilde{m}_f=6$) in the IMEX-E, IMEX-T and IMEX-E-Trap schems when solving the non-smooth case, and hence the CPU time of these schemes is slightly longer than that of the TS-I and TS-III schemes without correction terms with the same step size. However, to reach the same accuracy, the first three schemes are running faster than the TS-I and TS-III schemes for both non-smooth and smooth cases. In addition, neither the stiff problem \eqref{eq:fdemultiterm} with non-smooth solution nor with smooth solution can be solved by the predictor-corrector scheme. Figure \ref{fig:comp4_1} shows the asymptotic relative error of all schemes for both non-smooth case and smooth case, which illustrates that the IMEX-T, IMEX-E and IMEX-E-Trap schemes are superior for solving the multi-term nonlinear stiff problem \eqref{eq:fdemultiterm} with smooth and non-smooth solutions, compared to the TS-I scheme, TS-III scheme and the predictor-corrector scheme. {We observe that numerical solutions produced by the predictor-corrector scheme blow up in this example, which confirms that the IMEX schemes we present have better stability compared to the explicit schemes. For more information on the linear stability of the predictor-corrector scheme, see \cite{Gar10}.}

\begin{table}[!h]
	\centering
	\caption{
		A comparison of relative errors and convergence rates of the IMEX-T, IMEX-E, IMEX-E-Trap,TS-I, TS-III and PC schemes for multiterm nonlinear stiff system \eqref{eq:fdemultiterm} with non-smooth solution (Case I) and smooth solution (Case II) (Example \ref{exm:3}), $T=1.$}\label{tbl:nlinear-comp4}
	\scalebox{0.85}{
		\begin{tabular}{c|c|ccc|ccc}
			\hline
			\hline
			\multirow{2}{*}{scheme}& \multirow{2}{*}{$h$}& \multicolumn{3}{|c|}{Non-smooth solution (Case I)} & \multicolumn{3}{|c}{Smooth solution (Case II)}\\
			
			\cline{3-5}
			\cline{6-8}

			~&~& $ E^r_\infty(h)$ & order & cpu time &  $ E^r_\infty(h)$ & order  & cpu time  \\
			\hline
			IMEX-T& $2^{-5}$  &	1.30e-04	&	2.56	&	0.05 & 3.84e-04	&	2.31	&	0.03\\
			& $2^{-6}$  &	2.19e-05	&	2.75	&	0.07 & 7.73e-05	&	2.27	&	0.06\\
			& $2^{-7}$  &	3.25e-06	&	2.11	&	0.13 & 1.61e-05	&	2.20	&	0.11\\
			& $2^{-8}$  &	7.51e-07	&	*	&	0.31 & 3.49e-06	&	*	&	0.26\\
			\hline
			IMEX-E& $2^{-5}$  &	2.05e-04	&	2.40	&	0.06 & 7.52e-04 & 2.40	&	0.04\\
			& $2^{-6}$  &	3.89e-05	&	2.39	&	0.08 & 1.42e-04 & 2.40	&	0.07\\
			& $2^{-7}$  &	7.44e-06	&	2.35	&	0.16 & 2.69e-05 & 2.38	&	0.14\\
			& $2^{-8}$  &	1.46e-06	&	 *       &	0.32 & 5.18e-06 & 	*	&	0.30\\
			\hline
   IMEX-E-Trap& $2^{-5}$ &	2.11e-04	&	2.40	&	0.05 & 8.36e-04	&	2.38	&	0.03 \\
			& $2^{-6}$ &	4.02e-05	&	2.37	&	0.07 & 1.61e-04	&	2.37	&	0.06\\
			& $2^{-7}$ &	7.79e-06	&	2.32	&	0.16 & 3.10e-05	&	2.35	&	0.14\\
			& $2^{-8}$ &	1.56e-06	&	*	&	0.29     & 6.10e-06	&	*	&	0.26\\
			\hline
			TS-I& $2^{-5}$&	1.10e-01	&	0.57	&	0.05& 1.42e-02	&	1.49	&	0.05\\
			& $2^{-6}$&	7.39e-02	&	0.57	&	0.06& 5.05e-03	&	1.49	&	0.05\\
			& $2^{-7}$&	4.97e-02	&	0.57	&	0.13& 1.80e-03	&	1.48	&	0.11\\
			& $2^{-8}$&	3.35e-02	&	*	&	0.25& 6.46e-04	&	*	&	0.23\\
			\hline
			TS-III& $2^{-5}$&	1.10e-03	&	0.67	&	0.01& 3.27e-03	&	1.54	&	0.01\\
			& $2^{-6}$&	6.93e-04	&	0.68	&	0.04& 1.12e-03	&	1.52	&	0.03\\
			& $2^{-7}$&	4.31e-04	&	0.54	&	0.05& 3.91e-04	&	1.51	&	0.04\\
			& $2^{-8}$&	2.97e-04	&	*	&	0.07& 1.37e-04	&	*	&	0.08\\
			\hline
		PC  & $2^{-5} $&	Calculation& Failed	&		&	 Calculation & Failed	&			\\
		 	& $2^{-6}$&	Calculation& Failed	&		&	 Calculation & Failed	&			\\
			& $2^{-7}$&	Calculation& Failed	&		&	 Calculation & Failed	&			\\
			& $2^{-8}$&	Calculation& Failed	&		&	 Calculation & Failed	&			\\
			\hline
		\end{tabular}}
	\end{table}

\section{Proofs}\label{sec:proof}
In this section, we provide the convergence analysis of the two methods  \eqref{scheme-ts-e} and \eqref{scheme-ts-t}. We first introduce two lemmas.

\begin{lem}[\cite{Lub86,ZengLLT13}]\label{lem:weights}
	Let $\{\omega^{(\beta)}_n\}$ be given by \eqref{eq:genfun}. Then we have $\omega^{(\beta)}_n>0,$ and $\omega^{(\beta)}_n=O(n^{\beta-1}),n>0$.	
\end{lem}

\begin{lem}\label{lem:weight2}
Let $m_u$ and $m_f$ be positive integers and
%and  $W_{n,k}^{(\beta,\sigma)}$ and $W_{n,k}^{(\beta,\delta)}$ be defined by \eqref{eq:van1}.
the discrete operators $I_{h,\sigma}^{\beta,n,m_u}$ and    $I_{h,\delta}^{\beta,n,m_f}$ be defined by \eqref{Ih}.
Suppose that $\{\sigma_r\}$ and $\{\delta_r\}$ are sequences of strictly increasing positive numbers. Then there exists a positive constant $C$ independent of $n$ such that
\begin{equation}\label{eq:weight1}
|W_{n,k}^{(\beta,\sigma)}|\leq C(n^{\sigma_{m_u}+\beta-2}+n^{\beta-1}),
\;\quad	|W_{n,k}^{(\beta,\delta)}|\leq C(n^{\delta_{m_f}+\beta-2}+n^{\beta-1}).
\end{equation}
\end{lem}

\begin{proof*}
Letting $\nu=\sigma_r$ in \eqref{eq:lubich-1}, we  derive
\begin{equation*}
\frac{\Gamma(\sigma_r+1)}{\Gamma(\sigma_r+1+\beta)}t_n^{\sigma_r+\beta}
-h^\beta \sum_{k=0}^n\omega_{n-k}^{(\beta)} t_k^{\sigma_r}
=O(h^2t_n^{\sigma_r+\beta-2})+O(h^{1+\sigma_r}t_n^{\beta-1}),
\end{equation*}
which is equivalent to
\begin{equation*}
\frac{\Gamma(\sigma_r+1)}{\Gamma(\sigma_r+1+\beta)}n^{\sigma_r+\beta}
- \sum_{k=0}^n\omega_{n-k}^{(\beta)} k^{\sigma_r}
=O(n^{\sigma_r+\beta-2})+O(n^{\beta-1}).
\end{equation*}
Applying \eqref{eq:van1} yields
\begin{equation*}
\sum_{k=1}^{m_{u}}W_{n,k}^{(\beta,\sigma)} k^{\sigma_r}
=O(n^{\sigma_r+\beta-2})+O(n^{\beta-1}),\quad r=1,2,...,m_u,
\end{equation*}
which leads to $W_{n,k}^{(\beta,\sigma)} = O\left(n^{\sigma_1+\beta-2}\right)+...+O\left(n^{\sigma_{m_u}+\beta-2}\right)
+O\left(n^{\beta-1}\right)$. Since  $\sigma_r<\sigma_{r+1}$, we have
$|W_{n,k}^{(\beta,\sigma)}|\leq C(n^{\sigma_{m_u}+\beta-2}+n^{\beta-1})$.
We can similarly obtain $|W_{n,k}^{(\beta,\delta)}|\leq C(n^{\delta_{m_f}+\beta-2}+n^{\beta-1})$, which ends the proof.	
\end{proof*}

%\begin{rem}
For $\widehat{W}_{n,k}^{(f)}$ defined by \eqref{eq:extraweights}, we can easily derive
\begin{equation*}
\sum_{k=1}^{\widetilde{m}_f}\widehat{W}^{(f)}_{n,k} k^{\delta_r}=n^{\delta_r}-2({n-1})^{\delta_r}+({n-2})^{\delta_r}=O(n^{\delta_r-2}),
\quad r=1,\cdots,\widetilde{m}_f,
\end{equation*}
which leads to
\begin{equation}\label{eq:wf-2}
|\widehat{W}^{(f)}_{n,k} |\leq Cn^{\delta_{\widetilde{m}_f}-2}.
\end{equation}

We can similarly derive
\begin{equation}\label{eq:weightf}
|\widetilde{W}^{(f)}_{n,k} |\leq Cn^{\delta_{\widetilde{m}_f}-2},{\qquad}
|\widetilde{W}^{(u)}_{n,k} |\leq Cn^{\sigma_{\widetilde{m}_u}-2},
\end{equation}
where $\widetilde{W}^{(f)}_{n,k}$ and $\widetilde{W}^{(u)}_{n,k}$ are defined by
\eqref{eq:zeng1-2} and  \eqref{eq:zeng2-2}, respectively.
%\end{rem}

\begin{lem}[\cite{Dixon85}]\label{gronwall}
Let  $0<\beta<1$ and $c_{n,k}$ satisfy $0\leq c_{n,k} \leq (n-k)^{\beta-1},\;0\leq k\leq n$.
Assume that $e_n$ satisfies
$$|e_n|\leq A_0+Mh^{\beta}\sum_{k=0}^{n-1}c_{n,k}|e_k|,{\quad}n=1,2,...,N,$$
where $A_0,M>0$ and $nh\leq T,\,T>0$. Then there exists a positive constant $C$ independent of $n,h$ such that
$$|e_n|\leq CA_0,{\quad}n=1,2,...,N.$$
\end{lem}
Next, we present the proof for Theorem \ref{thm:convergence}. \begin{proof*}
Let $e_n=u(t_n)-U_n$. Then from \eqref{eq:intf1} and  \eqref{scheme-ts-e}, we have the following error equation
\begin{eqnarray}\label{error-1}
e_n&=&\lambda  h^{\beta}\left[\sum_{k=0}^n\omega_{n-k}^{(\beta)}e_k
+\sum_{k=1}^{m_u} W_{n,k}^{(\beta,\sigma)} e_k\right] \notag\\
&&+h^{\beta}\left[\sum_{k=0}^{n-1}\omega_{n-k}^{(\beta)}(f(t_k,u(t_k))-f(t_k,U_k))
+\sum_{k=1}^{m_f} W_{n,k}^{(\beta,\delta)} (f(t_k,u(t_k))-f(t_k,U_k))\right] \notag \\
&&+h^{\beta}\omega_0^{(\beta)}\bigg[2\big(f(t_{n-1},u(t_{n-1}))-f(t_{n-1},U_{n-1})\big)
-\big(f(t_{n-2},u(t_{n-2}))-f(t_{n-2},U_{n-2})\big) \notag \\
&&+\sum_{k=1}^{\widetilde{m}_f} \widehat{W}_{n,k}^{(f)} (f(t_k,u(t_k))-f(t_k,U_k))\bigg]+ R^n_E,
\end{eqnarray}
where $R^n_E$ is defined in  \eqref{eq:intf1}.
By simple calculation, we can derive
$$|R^n_E|\leq Ch^q,{\quad}
q={\min\{2,\sigma_{m_u+1}+\beta,\delta_{m_f+1}+\beta,\delta_{\widetilde{m}_f+1}+\beta\}}.$$

Since $f(t,u)$ satisfies the Lipschitz condition with respect to the second argument $u$,
i.e., $|f(t,x)-f(t,y)|\leq L|x-y|$, $L>0$,
we have  from \eqref{error-1} that
\begin{eqnarray}\label{error-2}
|e_n| &\leq &|\lambda | h^{\beta}\left[\sum_{k=0}^n\omega_{n-k}^{(\beta)}|e_k|
+\sum_{k=1}^{m_u} |W_{n,k}^{(\beta,\sigma)}| |e_k|\right]
 +Lh^{\beta}\left[\sum_{k=0}^{n-1}\omega_{n-k}^{(\beta)}|e_k|
+\sum_{k=1}^{m_f} |W_{n,k}^{(\beta,\delta)}| |e_k|\right] \notag \\
&&+Lh^{\beta}\omega_0^{(\beta)}\bigg[2|e_{n-1}|+|e_{n-2}|
+\sum_{k=1}^{\widetilde{m}_f} |\widehat{W}_{n,k}^{(f)}| |e_k|\bigg]+ |R^n_E| \notag\\
&\leq&|\lambda| h^{\beta}\omega_{0}^{(\beta)}|e_n|+
h^{\beta} \sum_{k=0}^{n-1}c_{n,k}|e_k|
+h^{\beta}  \sum_{k=1}^{m}W_{n,k} |e_k| + |R^n_E|,
\end{eqnarray}
where $c_{n,k}=(|\lambda | +L)\omega_{n-k}^{(\beta)}\,(0\leq k \leq n-3)$,
$c_{n,n-2}=(|\lambda|  +L)\omega_{2}^{(\beta)} +L \omega_0^{(\beta)}$,
$c_{n,n-1}=(|\lambda|  +L)\omega_{1}^{(\beta)} +2L \omega_0^{(\beta)}$,
$m=\max\{{m}_u,{m}_f,\widetilde{m}_f\}$, and $W_{n,k}\geq 0$ satisfies
$$W_{n,k}=|\lambda|  |W_{n,k}^{(\beta,\sigma)}|
+L|W_{n,k}^{(\beta,\delta)}| +L\omega_0^{(\beta)}
 |\widehat{W}_{n,k}^{(f)}|\leq Cn^{\beta}$$
 when $\sigma_{m_u},\delta_{m_f}\leq 2$ and $\delta_{\widetilde{m}_f}\leq 2+\beta$,
 see Eqs. \eqref{eq:weight1} and \eqref{eq:wf-2}.

We rewrite \eqref{error-2} into the following form
\begin{equation}\label{error-3}\begin{aligned}
(1-|\lambda | \omega_{0}^{(\beta)}h^{\beta})|e_n|
\leq & h^{\beta} \sum_{k=0}^{n-1}c_{n,k}|e_k|
+h^{\beta}  \sum_{k=1}^{m}W_{n,k} |e_k| + |R^n_E|.
\end{aligned}
\end{equation}
Since $\omega_{n}^{(\beta)}=O(n^{\beta-1})$ (see Lemma \ref{lem:weights}), we always have
$c_{n,k}\leq C(n-k)^{\beta-1}$.
Applying the generalized Gronwall's inequality (see Lemma \ref{gronwall}), we derive
\begin{equation}\label{error-4}\begin{aligned}
|e_n|\leq C\left(h^{\beta}  \sum_{k=1}^{m}W_{n,k} |e_k| + |R^n_E|\right).
\end{aligned}
\end{equation}
Note that $|R^n_E|\leq Ch^q$
and $W_{n,k}\leq Cn^{\beta}$, we have $h^{\beta}W_{n,k}\leq Ct_n^{\beta}\leq CT^{\beta}$
that leads to
$|e_n|\leq C\left( \sum_{k=1}^{m}|e_k| + h^q\right)$, which completes the proof.
\end{proof*}

The proof of Theorem \ref{thm:convergence2} is very similar, and hence it is omitted here.

\section{Conclusion}
We proposed two second-order IMEX schemes (see IMEX-E of \eqref{scheme-ts-e} and IMEX-T of \eqref{scheme-ts-t})  for  nonlinear FODEs with non-smooth solutions by using suitable  correction terms. We proved the convergence and linear stability of the IMEX-E and  IMEX-T schemes. The stability region of the IMEX-E scheme is bounded,  while the IMEX-T scheme is $A(\frac{\beta\pi}{2})$-stable, that is, the IMEX-T scheme is unconditionally   stable.

In order to obtain the derived IMEX schemes,
we presented the strategies of utilizing  suitable correction terms both in the approximation of fractional integrals and in extrapolation or Taylor expansion which are adopted  to linearize the schemes. The correction terms are useful to keep second-order accuracy of the IMEX schemes for solving nonlinear/stiff FODEs with non-smooth solutions. We further considered  the extension of these strategies to construct other second-order schemes with correction terms for the considered FODEs. We also discussed how to extend the present IMEX schemes from single-term FODEs to multi-term FODEs and systems; see Section \ref{sec4}.

We provided numerical examples to verify the efficiency of the proposed schemes, which shows second-order convergence for both smooth and non-smooth solutions by choosing suitable correction terms when solving a stiff system, a nonlinear FODE and a stiff nonlinear multi-term fractional differential system. {It was observed that when solving problems with non-smooth solution, applying suitable correction terms can significantly improve the accuracy, however, excessive use of correction terms is not conductive to raising accuracy, especially for the long-term simulation. Moreover, comparison between the present schemes and the existing ones illustrated that for the same level of accuracy, the present schemes cost less computational time for both smooth and non-smooth solutions.

In future work, we will focus on boundary value problems of fractional differential equations with non-smooth solution and propose high-order numerical methods. }

%The strategies presented in this paper, of using extrapolation or Taylor expansion with correction terms to obtain the IMEX schemes  and applying correction terms to increase accuracy in numerical schemes for solving nonlinear TFDEs with non-smooth solution, can be also utilized in other existing schemes or constructing new schemes to obtain high accuracy. As a result, we construct the IMEX-E-Trap scheme \eqref{scheme-ts-trap}; see  Section 4.

%Moreover, the proposed schemes, including the IMEX-E, IMEX-T and IMEX-E-Trap schemes can be applied to stiff time-fractional differential systems and multi-term nonlinear time-fractional differential equations/systems.  Example \ref{exm:3} shows that for solving a stiff nonlinear multi-term time-fractional differential system,  the proposed schemes can be of second-order convergence for both smooth solutions and non-smooth solutions.

\section*{Acknowledgment}
This work was supported by the MURI/ARO on ``Fractional PDEs for Conservation Laws and Beyond: Theory, Numerics and  Applications  (W911NF-15-1-0562)'',  and  also  by  NSF (DMS  1216437). The first author was also partially supported by NSF of China (No.11271036) and the third author of this work was also partially supported by a start-up fund from WPI.

%%%%%%%%%%%%%%%%%%%%%%%%%%%%%%%%%%%%%%%%%%%

%%%%%%%%%%%%%%%%%%%%%%%%%%%%%%%%%%%%%%%%%%%
\end{document}